# MISSPECIFICATION IN INFINITE-DIMENSIONAL BAYESIAN STATISTICS

BY B. J. K. KLEIJN AND A. W. VAN DER VAART

*Vrije Universiteit Amsterdam*

We consider the asymptotic behavior of posterior distributions if the model is misspecified. Given a prior distribution and a random sample from a distribution $P_0$, which may not be in the support of the prior, we show that the posterior concentrates its mass near the points in the support of the prior that minimize the Kullback–Leibler divergence with respect to $P_0$. An entropy condition and a prior-mass condition determine the rate of convergence. The method is applied to several examples, with special interest for infinite-dimensional models. These include Gaussian mixtures, nonparametric regression and parametric models.

**1. Introduction.** Of all criteria for statistical estimation, asymptotic consistency is among the least disputed. Consistency requires that the estimation procedure come arbitrarily close to the true, underlying distribution, if enough observations are used. It is of a frequentist nature, because it presumes a notion of an underlying, true distribution for the observations. If applied to posterior distributions, it is also considered a useful property by many Bayesians, as it could warn one away from prior distributions with undesirable, or unexpected, consequences. Priors which lead to undesirable posteriors have been documented, in particular, for non- or semiparametric models (e.g., [4, 5]), in which case it is also difficult to motivate a particular prior on purely intuitive, subjective grounds.

In the present paper we consider the situation where the posterior distribution cannot possibly be asymptotically consistent, because the model, or the prior, is misspecified. From a frequentist point of view, the relevance of studying misspecification is clear, because the assumption that the model contains the true, underlying distribution may lack realistic motivation in many practical situations. From an objective Bayesian point of view, the question is of interest, because, in principle, the Bayesian paradigm allows unrestricted choice of a prior, and, hence, we must allow for the possibility that the fixed distribution of the observations does not belong to the support of the prior. In this paper we show that in such a case the posterior will concentrate near a point in the support of the prior that is closest to the true sampling distribution as measured through the Kullback–Leibler divergence, and we give a characterization for the rate of concentration near this point.

Throughout the paper we assume that $X_1, X_2, \ldots$ are i.i.d. observations, each distributed according to a probability measure $P_0$. Given a model $\mathscr{P}$ and a prior $\Pi$, supported on $\mathscr{P}$, the posterior mass of a measurable subset $B \subset \mathscr{P}$ is given by

$$(1.1) \quad \Pi_n(B|X_1, \ldots, X_n) = \int_B \prod_{i=1}^n p(X_i)\, d\Pi(P) \Big/ \int_{\mathscr{P}} \prod_{i=1}^n p(X_i)\, d\Pi(P).$$

Here it is assumed that the model is dominated by a $\sigma$-finite measure $\mu$, and the density of a typical element $P \in \mathscr{P}$ relative to the dominating measure is written $p$ and assumed appropriately measurable. If we assume that the model is well specified, that is, $P_0 \in \mathscr{P}$, then posterior consistency means that the posterior distributions concentrate an arbitrarily large fraction of their total mass in arbitrarily small neighborhoods of $P_0$, if the number of observations used to determine the posterior is large enough. To formalize this, we let $d$ be a metric on $\mathscr{P}$ and say that the Bayesian procedure for the specified prior is *consistent*, if, for every $\varepsilon > 0$, $\Pi_n(\{P : d(P, P_0) > \varepsilon\} | X_1, \ldots, X_n) \to 0$, in $P_0$-probability. More specific information concerning the asymptotic behavior of an estimator is given by its rate of convergence. Let $\varepsilon_n > 0$ be a sequence that decreases to zero and suppose that, for any constants $M_n \to \infty$,

$$(1.2) \qquad \Pi_n(P \in \mathscr{P} : d(P, P_0) > M_n \varepsilon_n | X_1, \ldots, X_n) \to 0,$$

in $P_0$-probability. The sequence $\varepsilon_n$ corresponds to a decreasing sequence of neighborhoods of $P_0$, the $d$-radius of which goes to zero with $n$, while still capturing most of the posterior mass. If (1.2) is satisfied, then we say that the rate of convergence is at least $\varepsilon_n$.

If $P_0$ is at a positive distance from the model $\mathscr{P}$ and the prior concentrates all its mass on $\mathscr{P}$, then the posterior is inconsistent as it will concentrate all its mass on $\mathscr{P}$ as well. However, in this paper we show that the posterior will still settle down near a given measure $P^* \in \mathscr{P}$, and we shall characterize the sequences $\varepsilon_n$ such that the preceding display is valid with $d(P, P^*)$ taking the place of $d(P, P_0)$.

One would expect the posterior to concentrate its mass near minimum Kullback–Leibler points, since asymptotically the likelihood $\prod_{i=1}^n p(X_i)$ is maximal near points of minimal Kullback–Leibler divergence. The integrand in the numerator of (1.1) is the likelihood, so subsets of the model in which the (log-)likelihood is large account for a large fraction of the total posterior mass. Hence, it is no great surprise that the appropriate point of convergence $P^*$ is a minimum Kullback–Leibler point in $\mathscr{P}$, but the general issue of rates (and which metric $d$ to use) turns out to be more complicated than expected. We follow the work by Ghosal, Ghosh and van der Vaart [8] for the well-specified situation, but need to adapt, change or extend many steps.

After deriving general results, we consider several examples in some detail, including Bayesian fitting of Gaussian mixtures using Dirichlet priors on the mixing distribution, the regression problem and parametric models. Our results on the regression problem allow one, for instance, to conclude that a Bayesian approach in the nonparametric problem that uses a prior on the regression function, but employs a normal distribution for the errors, will lead to consistent estimation of the regression function, even if the regression errors are non-Gaussian. This result, which is the Bayesian counterpart of the well-known fact that least squares estimators (the maximum likelihood estimators if the errors are Gaussian) perform well even if the errors are non-Gaussian, is important to validate the Bayesian approach to regression, but appears to have received little attention in the literature.

1.1. *Notation and organization.* Let $L_1(\mathscr{X}, \mathscr{A})$ denote the set of all finite signed measures on $(\mathscr{X}, \mathscr{A})$ and let $\mathrm{conv}(\mathscr{Q})$ be the convex hull of a set of measures $\mathscr{Q}$: the set of all finite linear combinations $\sum_i \lambda_i Q_i$ for $Q_i \in \mathscr{Q}$ and $\lambda_i \geq 0$ with $\sum_i \lambda_i = 1$. For a measurable function $f$, let $Qf$ denote the integral $\int f\, dQ$.

The paper is organized as follows. Section 2 contains the main results of the paper, in increasing generality. Sections 3, 4 and 5 concern the three classes of examples that we consider: mixtures, the regression model and parametric models. Sections 6 and 7 contain the proofs of the main results, where the necessary results on tests are developed in Section 6 and are of independent interest. The final section is a technical appendix.

**2. Main results.** Let $X_1, X_2, \ldots$ be an i.i.d. sample from a distribution $P_0$ on a measurable space $(\mathscr{X}, \mathscr{A})$. Given a collection $\mathscr{P}$ of probability distributions on $(\mathscr{X}, \mathscr{A})$ and a prior probability measure $\Pi$ on $\mathscr{P}$, the posterior measure is defined as in (1.1) (where $0/0 = 0$ by definition). Here it is assumed that the "model" $\mathscr{P}$ is dominated by a $\sigma$-finite measure $\mu$ and that $x \mapsto p(x)$ is a density of $P \in \mathscr{P}$ relative to $\mu$ such that the map $(x, p) \mapsto p(x)$ is measurable relative to the product of $\mathscr{A}$ and an appropriate $\sigma$-field on $\mathscr{P}$, so that the right-hand side of (1.1) is a measurable function of $(X_1, \ldots, X_n)$ and a probability measure as a function of $B$ for every $X_1, \ldots, X_n$ such that the denominator is positive. The "true" distribution $P_0$ may or may not belong to the model $\mathscr{P}$. For simplicity of notation, we assume that $P_0$ possesses a density $p_0$ relative to $\mu$ as well.

Informally we think of the model $\mathscr{P}$ as the "support" of the prior $\Pi$, but we shall not make this precise in a topological sense. At this point we only assume that the prior concentrates on $\mathscr{P}$ in the sense that $\Pi(\mathscr{P}) = 1$ (but we note later that this too can be relaxed). Further requirements are made in the statements of the main results. Our main theorems implicitly assume the existence of a point $P^* \in \mathscr{P}$ minimizing the Kullback–Leibler divergence of $P_0$ to the model $\mathscr{P}$. In particular, the minimal Kullback–Leibler divergence is assumed to be finite, that is, $P^*$ satisfies

$$(2.1) \qquad -P_0 \log \frac{p^*}{p_0} < \infty.$$

By the convention that $\log 0 = -\infty$, the above implies that $P_0 \ll P^*$ and, hence, we assume without loss of generality that the density $p^*$ is strictly positive at the observations.

Our theorems give sufficient conditions for the posterior distribution to concentrate in neighborhoods of $P^*$ at a rate that is determined by the amount of prior mass "close to" the minimal Kullback–Leibler point $P^*$ and the "entropy" of the model. To specify the terms between quotation marks, we make the following definitions.

We define the entropy and the neighborhoods in which the posterior is to concentrate its mass relative to a semi-metric $d$ on $\mathscr{P}$. The general results are formulated relative to an arbitrary semi-metric and next the conditions will be simplified for more specific choices. Whether or not these simplifications can be made depends on the model $\mathscr{P}$, convexity being an important special case (see Lemma 2.2). Unlike in the case of well-specified priors, considered, for example, in [8], the Hellinger distance is not always appropriate in the misspecified situation. The general entropy bound is formulated in terms of a *covering number for testing under misspecification*, defined for $\varepsilon > 0$ as follows: we define $N_t(\varepsilon, \mathscr{P}, d; P_0, P^*)$ as the minimal number $N$ of convex sets $B_1, \ldots, B_N$ of probability measures on $(\mathscr{X}, \mathscr{A})$ needed to cover the set $\{P \in \mathscr{P} : \varepsilon < d(P, P^*) < 2\varepsilon\}$ such that, for every $i$,

$$(2.2) \qquad \inf_{P \in B_i} \sup_{0 < \alpha < 1} -\log P_0 \left(\frac{p}{p^*}\right)^\alpha \geq \frac{\varepsilon^2}{4}.$$

If there is no finite covering of this type, we define the covering number to be infinite. We refer to the logarithms $\log N_t(\varepsilon, \mathscr{P}, d; P_0, P^*)$ as *entropy numbers for testing under misspecification*. These numbers differ from ordinary metric entropy numbers in that the covering sets $B_i$ are required to satisfy the preceding display rather than to be balls of radius $\varepsilon$. We insist that the sets $B_i$ be convex and that (2.2) hold for every $P \in B_i$. This implies that (2.2) may involve measures $P$ that do not belong to the model $\mathscr{P}$ if this is not convex itself.

For $\varepsilon > 0$, we define a specific kind of Kullback–Leibler neighborhood of $P^*$ by

$$(2.3) \quad B(\varepsilon, P^*; P_0) = \left\{ P \in \mathscr{P} : -P_0 \log \frac{p}{p^*} \leq \varepsilon^2, P_0 \left( \log \frac{p}{p^*} \right)^2 \leq \varepsilon^2 \right\}.$$

THEOREM 2.1. *For a given model $\mathscr{P}$, prior $\Pi$ on $\mathscr{P}$ and some $P^* \in \mathscr{P}$, assume that $-P_0 \log(p^*/p_0) < \infty$ and $P_0(p/p^*) < \infty$ for all $P \in \mathscr{P}$. Suppose that there exist a sequence of strictly positive numbers $\varepsilon_n$ with $\varepsilon_n \to 0$ and $n\varepsilon_n^2 \to \infty$ and a constant $L > 0$, such that, for all $n$,*

$$(2.4) \quad \Pi(B(\varepsilon_n, P^*; P_0)) \geq e^{-Ln\varepsilon_n^2},$$

$$(2.5) \quad N_t(\varepsilon, \mathscr{P}, d; P_0, P^*) \leq e^{n\varepsilon_n^2} \quad \text{for all } \varepsilon > \varepsilon_n.$$

*Then for every sufficiently large constant $M$, as $n \to \infty$,*

$$(2.6) \quad \Pi_n(P \in \mathscr{P} : d(P, P^*) \geq M\varepsilon_n | X_1, \ldots, X_n) \to 0 \quad \text{in } L_1(P_0^n).$$

The proof of this theorem is given in Section 7. The theorem does not explicitly require that $P^*$ be a point of minimal Kullback–Leibler divergence, but this is implied by the conditions (see Lemma 6.4 below). The theorem is extended to the case of nonunique minimal Kullback–Leibler points in Section 2.4.

The two main conditions of Theorem 2.1 are a prior mass condition (2.4) and an entropy condition (2.5), which can be compared to Schwartz' conditions for posterior consistency (see [13]), or the two main conditions for the well-specified situation in [8]. Below we discuss the background of these conditions in turn.

The prior mass condition (2.4) reduces to the corresponding condition for the correctly specified case in [8] if $P^* = P_0$. Because $-P_0 \log(p^*/p_0) < \infty$, we may rewrite the first inequality in the definition (2.3) of the set $B(\varepsilon, P^*; P_0)$ as

$$-P_0 \log \frac{p}{p_0} \leq -P_0 \log \frac{p^*}{p_0} + \varepsilon^2.$$

Therefore, the set $B(\varepsilon, P^*; P_0)$ contains only $P \in \mathscr{P}$ that are within $\varepsilon^2$ of the minimal Kullback–Leibler divergence with respect to $P_0$ over the model. The lower bound (2.4) on the prior mass of $B(\varepsilon, P^*; P_0)$ requires that the prior measure assign a certain minimal share of its total mass to Kullback–Leibler neighborhoods of $P^*$. As argued in [8], a rough understanding of the exact form of (2.4) for the "optimal" rate $\varepsilon_n$ is that an optimal prior spreads its mass "uniformly" over $\mathscr{P}$. In the proof of Theorem 2.1, the prior mass condition serves to lower-bound the denominator in the expression for the posterior.

The background of the entropy condition (2.5) is more involved, but can be compared to a corresponding condition in the well-specified situation given in Theorem 2.1 of [8]. The purpose of the entropy condition is to measure the complexity of the model, a larger entropy leading to a slower rate of convergence. The entropy used in [8] is either the ordinary metric entropy $\log N(\varepsilon, \mathscr{P}, d)$, or the local entropy $\log N(\varepsilon/2, \{P \in \mathscr{P} : \varepsilon < d(P, P_0) < 2\varepsilon\}, d)$. For $d$ the Hellinger distance, the minimal $\varepsilon_n$ satisfying $\log N(\varepsilon_n, \mathscr{P}, d) = n\varepsilon_n^2$ is roughly the fastest rate of convergence for estimating a density in the model $\mathscr{P}$ relative to $d$ obtainable by any method of estimation (cf. [2]). We are not aware of a concept of "optimal rate of convergence" if the model is misspecified, but a rough interpretation of (2.5) given (2.4) would be that in the misspecified situation the posterior concentrates near the closest Kullback–Leibler point at the optimal rate pertaining to the model $\mathscr{P}$.

Misspecification requires that the complexity of the model be measured in a different, somewhat complicated way. In examples, depending on the semi-metric $d$, the covering numbers $N_t(\varepsilon, \mathscr{P}, d; P_0, P^*)$ can be related to ordinary metric covering numbers $N(\varepsilon, \mathscr{P}, d)$. For instance, we show below (see Lemmas 2.1–2.3) that, if the model $\mathscr{P}$ is convex, then the numbers $N_t(\varepsilon, \mathscr{P}, d; P_0, P^*)$ are bounded by the covering numbers $N(\varepsilon, \mathscr{P}, d)$ if the distance $d(P_1, P_2)$ equals the Hellinger distance between the measures $Q_i$ defined by $dQ_i = (p_0/p^*) \, dP_i$, that is, a weighted Hellinger distance between $P_1$ and $P_2$.

In the well-specified situation we have $P^* = P_0$, and the entropy numbers for testing can be bounded above by ordinary entropy numbers for the Hellinger distance. Thus, Theorem 2.1 becomes a refinement of the main theorem of [8]. To see this, we first note that, since $-\log x \geq 1 - x$ for $x > 0$,

$$-\log P_0 \left(\frac{p}{p_0}\right)^{1/2} \geq 1 - \int \sqrt{p}\sqrt{p_0} \, d\mu = \frac{1}{2} h^2(p, p_0).$$

It follows that the left-hand side of (2.2) with $\alpha = 1/2$ and $p_0 = p^*$ is bounded below by $\inf_{P \in B_i} h^2(p, p_0)$. Because Hellinger balls are convex, they are eligible candidates for the sets $B_i$ required in the definition of the covering numbers for testing. If we cover the set $\{P \in \mathscr{P} : 2\varepsilon < h(P, P_0) < 4\varepsilon\}$ by a minimal set of Hellinger balls of radius $\varepsilon/2$, then these balls automatically satisfy (2.2), by the triangle inequality. It follows that

$$N_t(\varepsilon, \mathscr{P}, h; P_0, P_0) \leq N(\varepsilon/2, \{P \in \mathscr{P} : 2\varepsilon < h(P, P_0) < 4\varepsilon\}, h).$$

The right-hand side is exactly a local covering number of the type used by [8]. Because the entropy numbers for testing allow general convex sets rather than balls of a given diameter, they appear to be genuinely smaller in general than local covering numbers. (Notably, the convex sets need not satisfy a size restriction.) However, we do not know of any examples of gains in the setting we are interested in, so that in the well-specified case there appears to be no use for the extended covering numbers as defined by (2.2). In the general, misspecified situation they are essential, even for standard parametric models, such as the one-dimensional normal location model.

At a more technical level, the entropy condition of [8] ensures the existence of certain tests of the measures $P$ versus the true measure $P_0$. In the misspecified case it is necessary to compare the measures $P$ to the minimal Kullback–Leibler point $P^*$, rather than to $P_0$. It turns out that the appropriate comparison is not a test of the measures $P$ versus the measure $P^*$ in the ordinary sense of testing, but to test the measures $Q(P)$ defined by $dQ(P) = (p_0/p^*)\,dP$ versus the measure $P_0$ [see (7.4)]. With $\mathcal{Q}$ the set of measures $Q(P)$ where $P$ ranges over $\mathcal{P}$, this leads to consideration of minimax testing risks of the type

$$\inf_\phi \sup_{Q \in \mathcal{Q}} (P_0^n \phi + Q^n(1-\phi)),$$

where the infimum is taken over all measurable functions $\phi$ taking values in $[0,1]$. A difference with the usual results on minimax testing risks is that the measures $Q$ may not be probability measures (and may in fact be infinite in general).

Extending arguments of Le Cam and Birgé, we show in Section 6 that, for a convex set $\mathcal{Q}$, the minimax testing risk in the preceding display is bounded above by

$$(2.7) \qquad \inf_{0<\alpha<1} \sup_{Q \in \mathcal{Q}} \rho_\alpha(P_0, Q)^n,$$

where the function $\alpha \mapsto \rho_\alpha(P_0, Q)$ is the Hellinger transform $\rho_\alpha(P,Q) = \int p^\alpha q^{1-\alpha}\,d\mu$. For $Q = Q(P)$, the Hellinger transform reduces to the map

$$\alpha \mapsto \rho_{1-\alpha}(Q(P), P_0) = P_0(p/p^*)^\alpha,$$

also encountered in (2.2). If the inequality in (2.2) is satisfied, then $P_0(p/p^*)^\alpha \le e^{-\varepsilon^2/4}$ and, hence, the set of measures $Q(P)$ with $P$ ranging over $B_i$ can be tested with error probabilities bounded by $e^{-n\varepsilon^2/4}$. For $\varepsilon$ bounded away from zero, or converging slowly to zero, these probabilities are exponentially small, ensuring that the posterior does not concentrate on the "unlikely alternatives" $B_i$.

The testing bound (2.7) is valid for convex alternatives $\mathcal{Q}$, but the alternatives of interest $\{P \in \mathcal{P} : d(P, P^*) > M\varepsilon\}$ are complements of balls and, hence, typically not convex. A test function for nonconvex alternatives can be constructed using a covering of $\mathcal{P}$ by convex sets. The entropy condition (2.5) controls the size of this cover and, hence, the rate of convergence in misspecified situations is determined by the covering numbers $N_t(\varepsilon, \mathscr{P}, d; P_0, P^*)$. Because the validity of the theorem only relies on the existence of suitable tests, the entropy condition (2.5) could be replaced by a testing condition. To be precise, condition (2.5) can be replaced by the condition that the conclusion of Theorem 6.3 is satisfied with $D(\varepsilon) = e^{n\varepsilon_n^2}$.

2.1. *Distances and testing entropy.* Because the entropies for testing are somewhat abstract, it is useful to relate them to ordinary entropy numbers. For our examples, the bound given by the following lemma is useful. We assume that, for some fixed constants $c, C > 0$ and for every $m \in \mathbb{N}$, $\lambda_1, \ldots, \lambda_m \geq 0$ with $\sum_i \lambda_i = 1$ and every $P, P_1, \ldots, P_m \in \mathscr{P}$ with $d(P, P_i) \leq cd(P, P^*)$ for all $i$,

$$(2.8) \quad \sum_i \lambda_i d^2(P_i, P^*) - C \sum_i \lambda_i d^2(P_i, P) \leq \sup_{0 < \alpha < 1} - \log P_0 \left( \frac{\sum_i \lambda_i p_i}{p^*} \right)^\alpha.$$

LEMMA 2.1. *If* (2.8) *holds, then there exists a constant* $A > 0$ *depending only on* $c$ *and* $C$ *such that, for all* $\varepsilon > 0$, $N_t(\varepsilon, \mathscr{P}, d; P_0, P^*) \leq N(A\varepsilon, \{P \in \mathscr{P} : \varepsilon < d(P, P^*) < 2\varepsilon\}, d)$. *[Any constant* $A \leq (1/8) \wedge (1/4\sqrt{C}) \wedge (\frac{1}{2}c)$ *works.]*

PROOF. For a given constant $A > 0$, we can cover the set $\mathscr{P}_\varepsilon := \{P \in \mathscr{P} : \varepsilon < d(P, P^*) < 2\varepsilon\}$ with $N = N(A\varepsilon, \mathscr{P}_\varepsilon, d)$ balls of radius $A\varepsilon$. If the centers of these balls are not contained in $\mathscr{P}_\varepsilon$, then we can replace these $N$ balls by $N$ balls of radius $2A\varepsilon$ with centers in $\mathscr{P}_\varepsilon$ whose union also covers the set $\mathscr{P}_\varepsilon$. It suffices to show that (2.2) is valid for $B_i$ equal to the convex hull of a typical ball $B$ in this cover. Choose $2A < c$. If $P \in \mathscr{P}_\varepsilon$ is the center of $B$ and $P_i \in B$ for every $i$, then $d(P_i, P^*) \geq d(P, P^*) - 2A\varepsilon$ by the triangle inequality and, hence, by assumption (2.8), the left-hand side of (2.2) with $B_i = \text{conv}(B)$ is bounded below by $\sum_i \lambda_i ((\varepsilon - 2A\varepsilon)^2 - C(2A\varepsilon)^2)$. This is bounded below by $\varepsilon^2/4$ for sufficiently small $A$. □

The logarithms $\log N(A\varepsilon, \{P \in \mathscr{P} : \varepsilon < d(P, P^*) < 2\varepsilon\}, d)$ of the covering numbers in the preceding lemma are called "local entropy numbers" and also the *Le Cam dimension* of the model $\mathscr{P}$ relative to the semi-metric $d$. They are bounded above by the simpler ordinary entropy numbers $\log N(A\varepsilon, \mathscr{P}, d)$. The preceding lemma shows that the entropy condition (2.5) can be replaced by the ordinary entropy condition $\log N(\varepsilon_n, \mathscr{P}, d) \leq n\varepsilon_n^2$ whenever the semi-metric $d$ satisfies (2.8).

If we evaluate (2.8) with $m = 1$ and $P_1 = P$, then we obtain, for every $P \in \mathscr{P}$,

$$(2.9) \quad d^2(P, P^*) \leq \sup_{0 < \alpha < 1} - \log P_0 \left( \frac{p}{p^*} \right)^\alpha.$$

(Up to a factor 16, this inequality is also implied by finiteness of the covering numbers for testing.) This simpler condition gives an indication about the metrics $d$ that may be used in combination with ordinary entropy. In Lemma 2.2 we show that if $d$ and the model $\mathscr{P}$ are convex, then the simpler condition (2.9) is equivalent to (2.8).

Because $-\log x \geq 1 - x$ for every $x > 0$, we can further simplify by bounding minus the logarithm in the right-hand side by $1 - P_0(p/p^*)^\alpha$. This yields the bound

$$d^2(P, P^*) \leq \sup_{0 < \alpha < 1} \left[ 1 - P_0 \left( \frac{p}{p^*} \right)^\alpha \right].$$

In the well-specified situation we have $P_0 = P^*$ and the right-hand side for $\alpha = 1/2$ becomes $1 - \int \sqrt{p}\sqrt{p_0}\, d\mu$, which is $1/2$ times the Hellinger distance between $P$ and $P_0$. In misspecified situations this method of lower bounding can be useless, as $1 - P_0(p/p^*)^\alpha$ may be negative for $\alpha = 1/2$. On the other hand, a small value of $\alpha$ may be appropriate, as it can be shown that as $\alpha \downarrow 0$ the expression $1 - P_0(p/p^*)^\alpha$ is proportional to the difference of Kullback–Leibler divergences $P_0 \log(p^*/p)$, which is positive by the definition of $P^*$. If this approximation can be made uniform in $p$, then a semi-metric $d$ which is bounded above by the Kullback–Leibler divergence can be used in the main theorem. We discuss this further in Section 6 and use this in the examples of Sections 4 and 5.

The case of convex models $\mathscr{P}$ is of interest, in particular, for non- or semi-parametric models, and permits some simplification. For a convex model, the point of minimal Kullback–Leibler divergence (if it exists) is automatically unique (up to redefinition on a null-set of $P_0$). Moreover, the expectations $P_0(p/p^*)$ are automatically finite, as required in Theorem 2.1, and condition (2.8) is satisfied for a weighted Hellinger metric. We show this in Lemma 2.3, after first showing that validity of the simpler lower bound (2.9) on the convex hull of $\mathscr{P}$ (if the semi-metric $d$ is defined on this convex hull) implies the bound (2.8).

LEMMA 2.2. *If $d$ is defined on the convex hull of $\mathscr{P}$, the maps $P \mapsto d^2(P, P')$ are convex on $\mathrm{conv}(\mathscr{P})$ for every $P' \in \mathscr{P}$ and (2.9) is valid for every $P$ in the convex hull of $\mathscr{P}$, then (2.8) is satisfied for $\frac{1}{2}d$ instead of $d$.*

LEMMA 2.3. *If $\mathscr{P}$ is convex and $P^* \in \mathscr{P}$ is a point at minimal Kullback–Leibler divergence with respect to $P_0$, then $P_0(p/p^*) \leq 1$ for every $P \in \mathscr{P}$ and (2.8) is satisfied with*

$$d^2(P_1, P_2) = \frac{1}{4} \int (\sqrt{p_1} - \sqrt{p_2})^2 \frac{p_0}{p^*}\, d\mu.$$

PROOF OF LEMMA 2.2. For the proof of Lemma 2.2, we first apply the triangle inequality repeatedly to find

$$\sum_i \lambda_i d^2(P_i, P^*) \le 2 \sum_i \lambda_i d^2(P_i, P) + 2 d^2(P, P^*)$$

$$\le 2 \sum_i \lambda_i d^2(P_i, P) + 4 d^2\left(P, \sum_i \lambda_i P_i\right) + 4 d^2\left(\sum_i \lambda_i P_i, P^*\right)$$

$$\le 6 \sum_i \lambda_i d^2(P_i, P) + 4 d^2\left(\sum_i \lambda_i P_i, P^*\right),$$

by the convexity of $d^2$. It follows that $d^2(\sum_i \lambda_i P_i, P^*) \ge (1/4) \sum_i \lambda_i d^2(P_i, P^*) - 3/2 \sum_i \lambda_i d^2(P_i, P)$. If (2.9) holds for $P = \sum_i \lambda_i P_i$, then we obtain (2.8) with $d^2$ replaced by $d^2/4$ and $C = 6$. □

PROOF OF LEMMA 2.3. For $P \in \mathscr{P}$, define a family of convex combinations $\{P_\lambda : \lambda \in [0,1]\} \subset \mathscr{P}$ by $P_\lambda = \lambda P + (1-\lambda) P^*$. For all values of $\lambda \in [0,1]$,

$$(2.10) \qquad 0 \le f(\lambda) := -P_0 \log \frac{p_\lambda}{p^*} = -P_0 \log\left(1 + \lambda\left(\frac{p}{p^*} - 1\right)\right),$$

since $P^* \in \mathscr{P}$ is at minimal Kullback–Leibler divergence with respect to $P_0$ in $\mathscr{P}$ by assumption. For every fixed $y \ge 0$, the function $\lambda \mapsto \log(1+\lambda y)/\lambda$ is nonnegative and increases monotonically to $y$ as $\lambda \downarrow 0$. The function is bounded in absolute value by 2 for $y \in [-1, 0]$ and $\lambda \le \frac{1}{2}$. Therefore, by the monotone and dominated convergence theorems applied to the positive and negative parts of the integrand in the right-hand side of (2.10),

$$f'(0+) = 1 - P_0\left(\frac{p}{p^*}\right).$$

Combining the fact that $f(0) = 0$ with (2.10), we see that $f'(0+) \ge 0$ and, hence, we find $P_0(p/p^*) \le 1$. The first assertion of Lemma 2.3 now follows.

For the proof that (2.9) is satisfied, we first note that $-\log x \ge 1 - x$, so that it suffices to show that $1 - P_0(p/p^*)^{1/2} \ge d^2(P, P^*)$. Now

$$\int (\sqrt{p^*} - \sqrt{p})^2 \frac{p_0}{p^*} d\mu = 1 + P_0 \frac{p}{p^*} - 2 P_0 \sqrt{\frac{p}{p^*}} \le 2 - 2 P_0 \sqrt{\frac{p}{p^*}},$$

by the first part of the proof. □

2.2. *Extensions*. In this section we give some generalizations of Theorem 2.1. Theorem 2.2 enables us to prove that optimal rates are achieved in parametric models. Theorem 2.3 extends Theorem 2.1 to situations in which the model, the prior and the point $P^*$ are dependent on $n$. Third, we consider the case in which the priors $\Pi_n$ assign a mass slightly less than 1 to the models $\mathscr{P}_n$.

Theorem 2.1 does not give the optimal rate of convergence $1/\sqrt{n}$ for finite-dimensional models $\mathscr{P}$, both because the choice $\varepsilon_n = 1/\sqrt{n}$ is excluded (by the condition $n\varepsilon_n^2 \to \infty$) and because the prior mass condition is too restrictive. The following theorem remedies this, but is more complicated. The adapted prior mass condition takes the following form: for all natural numbers $n$ and $j$,

$$(2.11) \qquad \frac{\Pi(P \in \mathscr{P} : j\varepsilon_n < d(P, P^*) < 2j\varepsilon_n)}{\Pi(B(\varepsilon_n, P^*; P_0))} \leq e^{n\varepsilon_n^2 j^2/8}.$$

THEOREM 2.2. *For a given model $\mathscr{P}$, prior $\Pi$ on $\mathscr{P}$ and some $P^* \in \mathscr{P}$, assume that $-P_0 \log(p^*/p_0) < \infty$ and $P_0(p/p^*) < \infty$ for all $P \in \mathscr{P}$. If $\varepsilon_n$ are strictly positive numbers with $\varepsilon_n \to 0$ and $\liminf n\varepsilon_n^2 > 0$, such that (2.11) and (2.5) are satisfied, then, for every sequence $M_n \to \infty$, as $n \to \infty$,*

$$(2.12) \quad \Pi_n(P \in \mathscr{P} : d(P, P^*) \geq M_n \varepsilon_n | X_1, \ldots, X_n) \to 0 \qquad \text{in } L_1(P_0).$$

There appears to be no compelling reason to choose the model $\mathscr{P}$ and the prior $\Pi$ the same for every $n$. The validity of the preceding theorems does not depend on this. We formalize this fact in the following theorem. For each $n$, we let $\mathscr{P}_n$ be a set of probability measures on $(\mathscr{X}, \mathscr{A})$ given by densities $p_n$ relative to a $\sigma$-finite measure $\mu_n$ on this space. Given a prior measure $\Pi_n$ on an appropriate $\sigma$-field, we define the posterior by (1.1) with $P^*$ and $\Pi$ replaced by $P_n^*$ and $\Pi_n$.

THEOREM 2.3. *The preceding theorems remain valid if $\mathscr{P}$, $\Pi$, $P^*$ and $d$ depend on $n$, but satisfy the given conditions for each $n$ (for a single constant $L$).*

As a final extension, we note that the assertion $P_0^n \Pi_n(P \in \mathscr{P}_n : d_n(P, P_n^*) \geq M_n \varepsilon_n | X_1, \ldots, X_n) \to 0$ of the preceding theorems remains valid even if the priors $\Pi_n$ do not put *all* their mass on the "models" $\mathscr{P}_n$ (but the models $\mathscr{P}_n$ do satisfy the entropy condition). Of course, in such cases the posterior puts mass outside the model and it is desirable to complement the above assertion with the assertion that $\Pi_n(\mathscr{P}_n | X_1, \ldots, X_n) \to 1$ in $L_1(P_0)$. The latter is certainly true if the priors put only very small fractions of their mass outside the models $\mathscr{P}_n$. More precisely, the latter assertion is true if

$$(2.13) \qquad \frac{1}{\Pi_n(B(\varepsilon_n, P_n^*, P_0))} \int_{\mathscr{P}_n^c} \left( P \frac{p_0}{p_n^*} \right)^n d\Pi_n(P) \leq o(e^{-2n\varepsilon_n^2}).$$

This observation is not relevant for the examples in the present paper. However, it may prove relevant to alleviate the entropy conditions in the preceding theorems in certain situations. These conditions limit the complexity of the models and it seems reasonable to allow a trade-off between complexity and prior mass. Condition (2.13) allows a crude form of such a trade-off: a small part $\mathscr{P}_n^c$ of the model may be more complex, provided that it receives a negligible amount of prior mass.

2.3. *Consistency.* The preceding theorems yield a rate of convergence $\varepsilon_n \to 0$, expressed as a function of prior mass and model entropy. In certain situations the prior mass and entropy may be hard to quantify. In contrast, for inferring consistency of the posterior, such quantification is unnecessary. This could be proved directly, as [13] achieved in the well-specified situation, but it can also be inferred from the preceding rate theorems. [A direct proof might actually give the same theorem with a slightly bigger set $B(\varepsilon, P^*; P_0)$.] We consider this for the situation of Theorem 2.1 only.

COROLLARY 2.1. *For a given model $\mathscr{P}$, prior $\Pi$ on $\mathscr{P}$ and some $P^* \in \mathscr{P}$, assume that $-P_0 \log(p^*/p_0) < \infty$ and $P_0(p/p^*) < \infty$ for all $P \in \mathscr{P}$. Suppose that, for every $\varepsilon > 0$,*

$$\Pi(B(\varepsilon, P^*; P_0)) > 0,$$

(2.14)

$$\sup_{\eta > \varepsilon} N_t(\eta, \mathscr{P}, d; P_0, P^*) < \infty.$$

*Then for every $\varepsilon > 0$, as $n \to \infty$,*

(2.15) $\quad \Pi_n(P \in \mathscr{P} : d(P, P^*) \geq \varepsilon | X_1, \ldots, X_n) \to 0 \qquad in \ L_1(P_0^n).$

PROOF. Define functions $f$ and $g$ by $f(\varepsilon) = \Pi(B(\varepsilon, P^*; P_0))$ and $g(\varepsilon) = \sup_{\eta > \varepsilon} N_t(\eta, \mathscr{P}, d; P_0, P^*)$. We shall show that there exists a sequence $\varepsilon_n \to 0$ such that $f(\varepsilon_n) \geq e^{-n\varepsilon_n^2}$ and $g(\varepsilon_n) \leq e^{n\varepsilon_n^2}$ for all sufficiently large $n$. This implies that the conditions of Theorem 2.1 are satisfied for this choice of $\varepsilon_n$ and, hence, the posterior converges with rate at least $\varepsilon_n$.

To show the existence of $\varepsilon_n$, define functions $h_n$ by

$$h_n(\varepsilon) = e^{-n\varepsilon^2}\left(g(\varepsilon) + \frac{1}{f(\varepsilon)}\right).$$

This is well defined and finite by the assumptions and $h_n(\varepsilon) \to 0$ as $n \to \infty$, for every fixed $\varepsilon > 0$. Therefore, there exists $\varepsilon_n \downarrow 0$ with $h_n(\varepsilon_n) \to 0$ [e.g., fix $n_1 < n_2 < \cdots$ such that $h_n(1/k) \leq 1/k$ for all $n \geq n_k$; next define $\varepsilon_n = 1/k$ for $n_k \leq n < n_{k+1}$]. In particular, there exists an $N$ such that $h_n(\varepsilon_n) \leq 1$ for $n \geq N$. This implies that $f(\varepsilon_n) \geq e^{-n\varepsilon_n^2}$ and $g(\varepsilon_n) \leq e^{n\varepsilon_n^2}$ for every $n \geq N$. □

2.4. *Multiple points of minimal Kullback–Leibler divergence.* In this section we extend Theorem 2.1 to the situation where there exists a set $\mathscr{P}^*$ of minimal Kullback–Leibler points.

Multiple minimal points can arise in two very different ways. First consider the situation where the true distribution $P_0$ and the elements of the model $\mathscr{P}$ possess different supports. Because the observations are sampled from $P_0$, they fall with probability one in the set where $p_0 > 0$ and, hence, the exact nature of the elements $p$ of the model $\mathscr{P}$ on the set $\{p_0 = 0\}$ is irrelevant. Clearly, multiple minima arise if the model contains multiple extensions of $P^*$ to the set on which $p_0 = 0$. In this case the observations do not provide the means to distinguish between these extensions and, consequently, no such distinction can be made by the posterior either. Theorems 2.1 and 2.2 may apply under this type of nonuniqueness, as long as we fix one of the minimal points, and the assertion of the theorem will be true for any of the minimal points as soon as it is true for one of the minimal points. This follows because, under the conditions of the theorem, $d(P_1^*, P_2^*) = 0$ whenever $P_1^*$ and $P_2^*$ agree on the set $p_0 > 0$, in view of (2.9).

Genuine multiple points of minimal Kullback–Leibler divergence may occur as well. For instance, one might fit a model consisting of normal distributions with means in $(-\infty, -1] \cup [1, \infty)$ and variance one, in a situation where the true distribution is normal with mean 0. The normal distributions with means $-1$ and $1$ both have the minimal Kullback–Leibler divergence. This situation is somewhat artificial and we are not aware of more interesting examples in the nonparametric or semiparametric case that interests us most in the present paper. Nevertheless, because it appears that the situation might arise, we give a brief discussion of an extension of Theorem 2.1.

Rather than to a single measure $P^* \in \mathscr{P}$, the extension refers to a finite subset $\mathscr{P}^* \subset \mathscr{P}$ of points at minimal Kullback–Leibler divergence. We give conditions under which the posterior distribution concentrates asymptotically near this set of points. We redefine our "covering numbers for testing under misspecification" $N_t(\varepsilon, \mathscr{P}, d; P_0, \mathscr{P}^*)$ as the minimal number $N$ of convex sets $B_1, \ldots, B_N$ of probability measures on $(\mathscr{X}, \mathscr{A})$ needed to cover the set $\{P \in \mathscr{P} : \varepsilon < d(P, \mathscr{P}^*) < 2\varepsilon\}$ such that

$$(2.16) \qquad \sup_{P^* \in \mathscr{P}^*} \inf_{P \in B_i} \sup_{0 < \alpha < 1} -\log P_0 \left(\frac{p}{p^*}\right)^\alpha \geq \frac{\varepsilon^2}{4}.$$

This roughly says that, for every $P \in \mathscr{P}$, there exists some minimal point to which we can apply arguments as before.

THEOREM 2.4. *For a given model $\mathscr{P}$, prior $\Pi$ on $\mathscr{P}$ and some subset $\mathscr{P}^* \subset \mathscr{P}$, assume that $-P_0 \log(p^*/p_0) < \infty$ and $P_0(p/p^*) < \infty$ for all $P \in \mathscr{P}$ and $P^* \in \mathscr{P}^*$. Suppose that there exist a sequence of strictly positive* numbers $\varepsilon_n$ with $\varepsilon_n \to 0$ and $n\varepsilon_n^2 \to \infty$ and a constant $L > 0$, such that, for all $n$ and all $\varepsilon > \varepsilon_n$,

$$\inf_{P^* \in \mathscr{P}^*} \Pi(B(\varepsilon_n, P^*; P_0)) \geq e^{-Ln\varepsilon_n^2}, \tag{2.17}$$

$$N_t(\varepsilon, \mathscr{P}, d; P_0, \mathscr{P}^*) \leq e^{n\varepsilon_n^2}. \tag{2.18}$$

Then for every sufficiently large constant $M > 0$, as $n \to \infty$,

$$\text{(2.19)} \quad \Pi_n(P \in \mathscr{P} : d(P, \mathscr{P}^*) \geq M\varepsilon_n | X_1, \ldots, X_n) \to 0 \qquad \text{in } L_1(P_0^n).$$

**3. Mixtures.** Let $\mu$ denote the Lebesgue measure on $\mathbb{R}$. For each $z \in \mathbb{R}$, let $x \mapsto p(x|z)$ be a fixed $\mu$-probability density on a measurable space $(\mathscr{X}, \mathscr{A})$ that depends measurably on $(x, z)$, and for a distribution $F$ on $\mathbb{R}$ define the $\mu$-density

$$p_F(x) = \int p(x|z) \, dF(z).$$

Let $P_F$ be the corresponding probability measure. In this section we consider mixture models $\mathscr{P} = \{P_F : F \in \mathscr{F}\}$, where $\mathscr{F}$ is the set of all probability distributions on a given compact interval $[-M, M]$. We consider consistency for general mixtures and derive a rate of convergence in the special case that the family $p(\cdot|z)$ is the normal location family, that is, with $\phi$ the standard normal density,

$$p_F(x) = \int \phi(x - z) \, dF(z). \tag{3.1}$$

The observations are an i.i.d. sample $X_1, \ldots, X_n$ drawn from a distribution $P_0$ on $(\mathscr{X}, \mathscr{A})$ with $\mu$-density $p_0$ which is not necessarily of the mixture form. As a prior for $F$, we use a Dirichlet process distribution (see [6, 7]) on $\mathscr{F}$.

3.1. *General mixtures.* We say that the model is $P_0$-identifiable if, for all pairs $F_1, F_2 \in \mathscr{F}$,

$$F_1 \neq F_2 \quad \implies \quad P_0(p_{F_1} \neq p_{F_2}) > 0.$$

Imposing this condition on the model excludes the first way in which nonuniqueness of $P^*$ may occur (as discussed in Section 2.4).

LEMMA 3.1. *Assume that* $-P_0 \log(p_F/p_0) < \infty$ *for some* $F \in \mathscr{F}$. *If the map* $z \mapsto p(x|z)$ *is continuous for all* $x$, *then there exists an* $F^* \in \mathscr{F}$ *that minimizes* $F \mapsto -P_0 \log(p_F/p_0)$ *over* $\mathscr{F}$. *If the model is* $P_0$-*identifiable, then this* $F^*$ *is unique.*

PROOF. If $F_n$ is a sequence in $\mathscr{F}$ with $F_n \to F$ for the weak topology on $\mathscr{F}$, then $p_{F_n}(x) \to p_F(x)$ for every $x$, since the kernel is continuous in $z$ (and, hence, also bounded as a result of the compactness of $[-M, M]$) and by use of the portmanteau lemma. Consequently, $p_{F_n} \to p_F$ in $L_1(\mu)$ by Scheffé's lemma. It follows that the map $F \mapsto p_F$ from $\mathscr{F}$ into $L_1(\mu)$ is continuous for the weak topology on $\mathscr{F}$. The set $\mathscr{F}$ is compact for this topology, by Prohorov's theorem. The Kullback–Leibler divergence $p \mapsto -P_0 \log(p/p_0)$ is lower semi-continuous as a map from $L_1(\mu)$ to $\mathbb{R}$ (see Lemma 8.2). Therefore, the map $F \mapsto -P_0 \log(p_F/p_0)$ is lower semi-continuous on the compactum $\mathscr{F}$ and, hence, attains its infimum on $\mathscr{F}$.

The map $p \mapsto -P_0 \log(p/p_0)$ is also convex. By the strict convexity of the function $x \mapsto -\log x$, we have, for any $\lambda \in (0,1)$,

$$-P_0 \log\left(\frac{\lambda p_1 + (1-\lambda)p_2}{p_0}\right) < -\lambda P_0 \log \frac{p_1}{p_0} - (1-\lambda) P_0 \log \frac{p_2}{p_0},$$

unless $P_0(p_1 = p_2) = 1$. This shows that the point of minimum of $F \mapsto P_0 \log(p_F/p_0)$ is unique if $\mathscr{F}$ is $P_0$-identifiable. $\square$

Thus, a minimal Kullback–Leibler point $P_{F^*}$ exists and is unique under mild conditions on the kernel $p$. Because the model is convex, Lemma 2.3 next shows that (2.9) is satisfied for the weighted Hellinger distance, whose square is equal to

$$(3.2) \qquad d^2(P_{F_1}, P_{F_2}) = \frac{1}{2} \int (\sqrt{p_{F_1}} - \sqrt{p_{F_2}})^2 \frac{p_0}{p_{F^*}} \, d\mu.$$

If $p_0/p_{F^*} \in L_\infty(\mu)$, then this expression is bounded by the squared Hellinger distance $H$ between the measures $P_{F_1}$ and $P_{F_2}$.

Because $\mathscr{F}$ is compact for the weak topology and the map $F \mapsto p_F$ from $\mathscr{F}$ to $L_1(\mu)$ is continuous (cf. the proof of Lemma 3.1), the model $\mathscr{P} = \{P_F : F \in \mathscr{F}\}$ is compact relative to the total variation distance. Because the Hellinger and total variation distances define the same uniform structure, the model is also compact relative to the Hellinger distance and, hence, it is totally bounded, that is, the covering numbers $N(\varepsilon, \mathscr{P}, H)$ are finite for all $\varepsilon$. Combined with the result of the preceding paragraph and Lemmas 2.2 and 2.3, this yields the result that the entropy condition of Corollary 2.1 is satisfied for $d$ as in (3.2) if $p_0/p_{F^*} \in L_\infty(\mu)$ and we obtain the following theorem.

THEOREM 3.1. *If $p_0/p_{F^*} \in L_\infty(\mu)$ and $\Pi(B(\varepsilon, P_{F^*}; P_0)) > 0$ for every $\varepsilon > 0$, then $\Pi_n(F : d(P_F, P_{F^*}) \geq \varepsilon | X_1, \ldots, X_n) \to 0$ in $L_1(P_0^n)$ for $d$ given by (3.2).*

3.2. *Gaussian mixtures.* Next we specialize to the situation where $p(x|z) = \phi(x - z)$ is a Gaussian convolution kernel and derive the rate of convergence. The Gaussian convolution model is well known to be $P_0$-identifiable if $P_0$ is Lebesgue absolutely continuous (see, e.g., [12]). Let $d$ be defined as in (3.2). We assume that $P_0$ is such that $-P_0 \log(p_F/p_0)$ is finite for some $F$, so that there exists a minimal Kullback–Leibler point $F^*$, by Lemma 3.1.

LEMMA 3.2. *If for some constant $C_1 > 0$, $d(p_{F_1}, p_{F_2}) \leq C_1 H(p_{F_1}, p_{F_2})$, then the entropy condition*

$$\log N(\varepsilon_n, \mathscr{P}, d) \leq n\varepsilon_n^2$$

*is satisfied for $\varepsilon_n$ a large enough multiple of $\log n/\sqrt{n}$.*

PROOF. Because the square of the Hellinger distance is bounded above by the $L_1$-norm, the assumption implies that $d^2(P_{F_1}, P_{F_2}) \leq C_1^2 \|P_{F_1} - P_{F_2}\|_1$. Hence, for all $\varepsilon > 0$, we have $N(C_1\varepsilon, \mathscr{P}, d) \leq N(\varepsilon^2, \mathscr{P}, \|\cdot\|_1)$. As a result of Lemma 3.3 in [9], there exists a constant $C_2 > 0$ such that

$$(3.3) \quad \|P_{F_1} - P_{F_2}\|_1 \leq C_2^2 \|P_{F_1} - P_{F_2}\|_\infty \max\left\{1, M, \sqrt{\log_+ \frac{1}{\|P_{F_1} - P_{F_2}\|_\infty}}\right\},$$

from which it follows that $N(C_2^2 \varepsilon \log(1/\varepsilon)^{1/2}, \mathscr{P}, \|\cdot\|_1) \leq N(\varepsilon, \mathscr{P}, \|\cdot\|_\infty)$ for small enough $\varepsilon$. With the help of Lemma 3.2 in [9], we see that there exists a constant $C_3 > 0$ such that, for all $0 < \varepsilon < e^{-1}$,

$$\log N(\varepsilon, \mathscr{P}, \|\cdot\|_\infty) \leq C_3 \left(\log \frac{1}{\varepsilon}\right)^2.$$

Combining all of the above, we note that, for small enough $\varepsilon > 0$,

$$\log N\left(C_1 C_2 \varepsilon^{1/2} \left(\log \frac{1}{\varepsilon}\right)^{1/4}, \mathscr{P}, d\right) \leq \log N(\varepsilon, \mathscr{P}, \|\cdot\|_\infty)$$

$$\leq C_3 \left(\log \frac{1}{\varepsilon}\right)^2.$$

So if we can find a sequence $\varepsilon_n$ such that, for all $n \geq 1$, there exists an $\varepsilon > 0$ such that

$$C_1 C_2 \varepsilon^{1/2} \left(\log \frac{1}{\varepsilon}\right)^{1/4} \leq \varepsilon_n \quad \text{and} \quad C_3 \left(\log \frac{1}{\varepsilon}\right)^2 \leq n\varepsilon_n^2,$$

then we have demonstrated that

$$\log N(\varepsilon_n, \mathscr{P}, d) \leq \log N\left(C_1 C_2 \varepsilon^{1/2} \left(\log \frac{1}{\varepsilon}\right)^{1/4}, \mathscr{P}, d\right)$$

$$\leq C_3 \left(\log \frac{1}{\varepsilon}\right)^2 \leq n\varepsilon_n^2.$$

One easily shows that this is the case for $\varepsilon_n = \max\{C_1 C_2, C_3\}(\log n/\sqrt{n})$ (in which case we choose, for fixed $n$, $\varepsilon = 1/n$), if $n$ is taken large enough. $\square$

We are now in position to apply Theorem 2.1. We consider, for given $M > 0$, the location mixtures (3.1) with the standard normal density $\phi$ as the kernel. We choose the prior $\Pi$ equal to a Dirichlet prior on $\mathscr{F}$ specified by a finite base measure $\alpha$ with compact support and positive, continuous Lebesgue density on $[-M, M]$.

THEOREM 3.2. *Let $P_0$ be a distribution on $\mathbb{R}$ dominated by Lebesgue measure $\mu$. Assume that $p_0/p_{F^*} \in L_\infty(\mu)$. Then the posterior distribution concentrates its mass around $P_{F^*}$ asymptotically at the rate $\log n/\sqrt{n}$ relative to the distance d on $\mathscr{P}$ given by (3.2).*

PROOF. The set of mixture densities $p_F$ with $F \in \mathscr{F}$ is bounded above and below by the upper and lower envelope functions

$$U(x) = \phi(x+M)\mathbb{1}_{\{x<-M\}} + \phi(x-M)\mathbb{1}_{\{x>M\}} + \phi(0)\mathbb{1}_{\{-M\leq x\leq M\}},$$
$$L(x) = \phi(x-M)\mathbb{1}_{\{x<0\}} + \phi(x+M)\mathbb{1}_{\{x\geq 0\}}.$$

So for any $F \in \mathscr{F}$,

$$P_0\left(\frac{p_F}{p_{F^*}}\right) \leq P_0 \frac{U}{L}$$
$$\leq \frac{\phi(0)}{\phi(2M)} P_0[-M, M]$$
$$\quad + P_0(e^{-2MX}\mathbb{1}_{\{X<-M\}} + e^{2MX}\mathbb{1}_{\{X>M\}}) < \infty,$$

because $p_0$ is essentially bounded by a multiple of $p_{F^*}$ and $P_{F^*}$ has sub-Gaussian tails. In view of Lemmas 2.2 and 2.3, the covering number for testing $N_t(\varepsilon, \mathscr{P}, d; P_0, P^*)$ in (2.5) is bounded above by the ordinary metric covering number $N(A\varepsilon, \mathscr{P}, d)$, for some constant $A$. Then Lemma 3.2 demonstrates that the entropy condition (2.5) is satisfied for $\varepsilon_n$ a large multiple of $\log n/\sqrt{n}$.

It suffices to verify the prior mass condition (2.4). Let $\varepsilon$ be given such that $0 < \varepsilon < e^{-1}$. By Lemma 3.2 in [9], there exists a discrete distribution function $F' \in D[-M, M]$ supported on at most $N \leq C_2 \log(1/\varepsilon)$ points $\{z_1, z_2, \ldots, z_N\} \subset [-M, M]$ such that $\|p_{F^*} - p_{F'}\|_\infty \leq C_1\varepsilon$, where $C_1, C_2 > 0$ are constants that depend on $M$ only. We write $F' = \sum_{j=1}^{N} p_j \delta_{z_j}$. Without loss of generality, we may assume that the set $\{z_j : j = 1, \ldots, N\}$ is $2\varepsilon$-separated. Namely, if this is not the case, we may choose a maximal $2\varepsilon$-separated subset of $\{z_j : j = 1, \ldots, N\}$ and shift the weights $p_j$ to the nearest point in the subset. A discrete $F''$ obtained in this fashion satisfies $\|p_{F'} - p_{F''}\|_\infty \leq 2\varepsilon\|\phi'\|_\infty$. So by virtue of the triangle inequality and the fact that the derivative of the standard normal kernel $\phi$ is bounded, a given $F'$ may be replaced by a $2\varepsilon$-separated $F''$ if the constant $C_1$ is changed accordingly.

By Lemma 3.3 in [9], there exists a constant $D_1$ such that the $L_1$-norm of the difference satisfies

$$\|P_{F^*} - P_{F'}\|_1 \leq D_1\varepsilon\left(\log\frac{1}{\varepsilon}\right)^{1/2},$$

for small enough $\varepsilon$. Using Lemma 3.6 in [9], we note, moreover, that there exists a constant $D_2$ such that, for any $F \in \mathscr{F}$,

$$\|P_F - P_{F'}\|_1 \leq D_2\left(\varepsilon + \sum_{j=1}^{N}|F[z_j - \varepsilon, z_j + \varepsilon] - p_j|\right).$$

So there exists a constant $D > 0$ such that, if $F$ satisfies $\sum_{j=1}^{N}|F[z_j - \varepsilon, z_j + \varepsilon] - p_j| \leq \varepsilon$, then

$$\|P_F - P_{F^*}\|_1 \leq D\varepsilon\left(\log\frac{1}{\varepsilon}\right)^{1/2}.$$

Let $Q(P)$ be the measure defined by $dQ(P) = (p_0/p_{F^*})\,dP$. The assumption that $p_0/p_{F^*}$ is essentially bounded implies that there exists a constant $K > 0$ such that $\|Q(P_{F_1}) - Q(P_{F_2})\|_1 \leq K\|P_{F_1} - P_{F_2}\|_1$ for all $F_1, F_2 \in \mathscr{F}$. Since $Q(P_{F^*}) = P_0$, it follows that there exists a constant $D' > 0$ such that, for small enough $\varepsilon > 0$,

$$\left\{F \in \mathscr{F} : \sum_{j=1}^{N}|F[z_j - \varepsilon, z_j + \varepsilon] - p_j| \leq \varepsilon\right\}$$
$$\subset \left\{F \in \mathscr{F} : \|Q(P_F) - P_0\|_1 \leq (D')^2\varepsilon\left(\log\frac{1}{\varepsilon}\right)^{1/2}\right\}.$$

We have that $dQ(P_F)/dP_0 = p_F/p_{F^*}$ and $P_0(p_F/p_{F^*}) \leq P_0(U/L) < \infty$. The Hellinger distance is bounded by the square root of the $L_1$-distance. Therefore, applying Lemma 8.1 with $\eta = \eta(\varepsilon) = D'\varepsilon^{1/2}(\log(1/\varepsilon))^{1/4}$, we see that the set of measures $P_F$ with $F$ in the set on the right-hand side of the last display is contained in the set

$$\left\{P_F : F \in \mathscr{F}, -P_0\log\frac{p_F}{p_{F^*}} \leq \zeta^2(\varepsilon), P_0\left(\log\frac{p_F}{p_{F^*}}\right)^2 \leq \zeta^2(\varepsilon)\right\}$$
$$\subset B(\zeta(\varepsilon), P_{F^*}; P_0),$$

where $\zeta(\varepsilon) = D''\eta(\varepsilon)(\log(1/\eta(\varepsilon))) \leq D''D'\varepsilon^{1/2}(\log(1/\varepsilon))^{5/4}$, for an appropriate constant $D''$, and small enough $\varepsilon$. It follows that

$$\Pi(B(\zeta(\varepsilon), P_{F^*}; P_0)) \geq \Pi\bigg\{F \in \mathscr{F} : \sum_{j=1}^{N} |F[z_j - \varepsilon, z_j + \varepsilon] - p_j| \leq \varepsilon\bigg\}.$$

Following [8] (Lemma 6.1) or Lemma A.2 in [9], we see that the prior measure at the right-hand side of the previous display is lower bounded by

$$c_1 \exp(-c_2 N \log(1/\varepsilon)) \geq \exp(-L(\log(1/\varepsilon))^2) \geq \exp(-L'(\log(1/\zeta(\varepsilon)))^2),$$

where $c_1 > 1$, $c_2 > 0$ are constants and $L = C_2 c_2 > 0$. So if we can find a sequence $\varepsilon_n$ such that, for each sufficiently large $n$, there exists an $\varepsilon > 0$ such that

$$\varepsilon_n \geq \zeta(\varepsilon), \qquad n\varepsilon_n^2 \geq \bigg(\log \frac{1}{\zeta(\varepsilon)}\bigg)^2,$$

then $\Pi(B(\varepsilon_n, P_{F^*}; P_0)) \geq \Pi(B(\zeta(\varepsilon), P_{F^*}; P_0)) \geq \exp(-L'n\varepsilon_n^2)$ and, hence, (2.4) is satisfied. One easily shows that, for $\varepsilon_n = \log n/\sqrt{n}$ and $\zeta(\varepsilon) = 1/\sqrt{n}$, the two requirements are fulfilled for sufficiently large $n$. $\square$

**4. Regression.** Let $P_0$ be the distribution of a random vector $(X, Y)$ satisfying $Y = f_0(X) + e_0$ for independent random variables $X$ and $e_0$ taking values in a measurable space $(\mathscr{X}, \mathscr{A})$ and in $\mathbb{R}$, respectively, and a measurable function $f_0 : \mathscr{X} \to \mathbb{R}$. The variables $X$ and $e_0$ have given marginal distributions, which may be unknown, but are fixed throughout the following. The purpose is to estimate the regression function $f_0$ based on a random sample of variables $(X_1, Y_1), \ldots, (X_n, Y_n)$ with the same distribution as $(X, Y)$.

A Bayesian approach to this problem might start from the specification of a prior distribution on a given class $\mathscr{F}$ of measurable functions $f : \mathscr{X} \to \mathbb{R}$. If the distributions of $X$ and $e_0$ are known, this is sufficient to determine a posterior. If these distributions are not known, then one might proceed to introduce priors for these unknowns as well. The approach we take here is to fix the distribution of $e_0$ to a normal or Laplace distribution, while aware of the fact that its true distribution may be different. We investigate the consequences of the resulting model misspecification. We shall show that misspecification of the error distribution does not have serious consequences for estimation of the regression function. In this sense a nonparametric Bayesian approach possesses the same robustness to misspecification as minimum contrast estimation using least squares or minimum absolute deviation. We shall also see that the use of the Laplace distribution requires no conditions on the tail of the distribution of the errors, whereas the normal distribution appears to give good results only if these tails are not too big. Thus, the tail robustness of minimum absolute deviation versus the nonrobustness of the method of least squares also extends to Bayesian regression.

We build the posterior based on a regression model $Y = f(X) + e$ for $X$ and $e$ independent, as is the assumption on the true distribution of $(X,Y)$. If we assume that the distribution $P_X$ of $X$ has a known form, then this distribution cancels out of the expression for the posterior on $f$. If, instead, we put independent priors on $f$ and $P_X$, respectively, then the prior on $P_X$ would disappear upon marginalization of the posterior of $(f, P_X)$ relative to $f$. Thus, for investigating the posterior for $f$, we may assume without loss of generality that the marginal distribution of $X$ is known. It can be absorbed into the dominating measure $\mu$ for the model.

For $f \in \mathscr{F}$, let $P_f$ be the distribution of the random variable $(X,Y)$ satisfying $Y = f(X) + e$ for $X$ and $e$ independent variables, $X$ having the same distribution as before and $e$ possessing a given density $p$, possibly different from the density of the true error $e_0$. We shall consider the cases that $p$ is normal and Laplace. Given a prior $\Pi$ on $\mathscr{F}$, the posterior distribution for $f$ is given by

$$B \mapsto \frac{\int_B \prod_{i=1}^n p(Y_i - f(X_i)) \, d\Pi(f)}{\int \prod_{i=1}^n p(Y_i - f(X_i)) \, d\Pi(f)}.$$

We shall show that this distribution concentrates near $f_0 + \mathrm{E}e_0$ in the case that $p$ is a normal density and near $f_0 + \mathrm{median}(e_0)$ if $p$ is Laplace, if these translates of the true regression function $f_0$ are contained in the model $\mathscr{F}$. If the prior is misspecified also in the sense that $f_0 + \mu \notin \mathscr{F}$ (where $\mu$ is the expectation or median of $e_0$), then, under some conditions, this remains true with $f_0$ replaced by a "projection" $f^*$ of $f_0$ on $\mathscr{F}$. In agreement with the notation in the rest of the paper, we shall denote the true distribution of an observation $(X,Y)$ by $P_0$ (stressing that, in general, $P_0$ is different from $P_f$ with $f = 0$). The model $\mathscr{P}$ as in the statement of the main results is the set of all distributions $P_f$ on $\mathscr{X} \times \mathbb{R}$ with $f \in \mathscr{F}$.

4.1. *Normal regression.* Suppose that the density $p$ is equal to the standard normal density $p(z) = (2\pi)^{-1/2} \exp(-\frac{1}{2}z^2)$. Then, with $\mu = \mathrm{E}e_0$,

(4.1)
$$\log \frac{p_f}{p_{f_0}}(X,Y) = -\frac{1}{2}(f - f_0)^2(X) + e_0(f - f_0)(X),$$

$$-P_0 \log \frac{p_f}{p_{f_0}} = \frac{1}{2}P_0(f - f_0 - \mu)^2 - \frac{1}{2}\mu^2.$$

It follows that the Kullback–Leibler divergence $f \mapsto -P_0 \log(p_f/p_0)$ is minimized for $f = f^* \in \mathscr{F}$ minimizing the map $f \mapsto P_0(f - f_0 - \mu)^2$.

In particular, if $f_0 + \mu \in \mathscr{F}$, then the minimizer is $f_0 + \mu$ and $P_{f_0 + \mu}$ is the point in the model that is closest to $P_0$ in the Kullback–Leibler sense. If also $\mu = 0$, then, even though the posterior on $P_f$ will concentrate asymptotically near $P_{f_0}$, which is typically not equal to $P_0$, the induced posterior on $f$ will concentrate near the true regression function $f_0$. This favorable property of Bayesian estimation is analogous to that of least squares estimators, also for nonnormal error distributions.

If $f_0 + \mu$ is not contained in the model, then the posterior for $f$ will in general not be consistent. We assume that there exists a unique $f^* \in \mathscr{F}$ that minimizes $f \mapsto P_0(f - f_0 - \mu)^2$, as is the case, for instance, if $\mathscr{F}$ is a closed, convex subset of $L_2(P_0)$. Under some conditions we shall show that the posterior concentrates asymptotically near $f^*$. If $\mu = 0$, then $f^*$ is the projection of $f_0$ into $\mathscr{F}$ and the posterior still behaves in a desirable manner. For simplicity of notation, we assume that $\mathrm{E}_0 e_0 = 0$.

The following lemma shows that (2.8) is satisfied for a multiple of the $L_2(P_0)$-distance on $\mathscr{F}$.

LEMMA 4.1. *Let $\mathscr{F}$ be a class of uniformly bounded functions $f : \mathscr{X} \to \mathbb{R}$ such that either $f_0 \in \mathscr{F}$ or $\mathscr{F}$ is convex and closed in $L_2(P_0)$. Assume that $f_0$ is uniformly bounded, that $\mathrm{E}_0 e_0 = 0$ and that $\mathrm{E}_0 e^{M|e_0|} < \infty$ for every $M > 0$. Then there exist positive constants $C_1, C_2, C_3$ such that, for all $m \in \mathbb{N}$, $f, f_1, \ldots, f_m \in \mathscr{F}$ and $\lambda_1, \ldots, \lambda_m \geq 0$ with $\sum_i \lambda_i = 1$,*

$$P_0 \log \frac{p_f}{p_{f^*}} \leq -\frac{1}{2} P_0 (f - f^*)^2,$$

(4.2) $$P_0 \left( \log \frac{p_{f^*}}{p_f} \right)^2 \leq C_1 P_0 (f - f^*)^2,$$

$$\sup_{0 < \alpha < 1} -\log P_0 \left( \frac{\sum_i \lambda_i p_{f_i}}{p_{f^*}} \right)^\alpha \geq C_2 \sum_i \lambda_i (P_0(f_i - f^*)^2 - C_3 P_0(f - f_i)^2).$$

PROOF. We have

(4.3) $$\log \frac{p_f}{p_{f^*}}(X, Y) = -\frac{1}{2}[(f_0 - f)^2 - (f_0 - f^*)^2](X) - e_0(f^* - f)(X).$$

The second term on the right-hand side has mean zero by assumption. The first term on the right-hand side has expectation $-\frac{1}{2} P_0 (f^* - f)^2$ if $f_0 = f^*$, as is the case if $f_0 \in \mathscr{F}$. Furthermore, if $\mathscr{F}$ is convex, the minimizing property of $f^*$ implies that $P_0(f_0 - f^*)(f^* - f) \geq 0$ for every $f \in \mathscr{F}$ and then the expectation of the first term on the right-hand side is bounded above by $-\frac{1}{2} P_0 (f^* - f)^2$. Therefore, in both cases (4.2) holds.

From (4.3) we also have, with $M$ a uniform upper bound on $\mathscr{F}$ and $f_0$,

$$P_0 \left( \log \frac{p_f}{p_{f^*}} \right)^2 \leq P_0 [(f^* - f)^2 (2f_0 - f - f^*)^2] + 2 P_0 e_0^2 P_0 (f^* - f)^2,$$

$$P_0\Big(\log\frac{p_f}{p_{f^*}}\Big)^2\Big(\frac{p_f}{p_{f^*}}\Big)^\alpha \le P_0[(f^*-f)^2(2f_0-f-f^*)^2+2e_0^2(f^*-f)^2]$$
$$\times e^{2\alpha(M^2+M|e_0|)}.$$

Both right-hand sides can be further bounded by a constant times $P_0(f-f^*)^2$, where the constant depends on $M$ and the distribution of $e_0$ only.

In view of Lemma 4.3 (below) with $p=p_{f^*}$ and $q_i=p_{f_i}$, we see that there exists a constant $C>0$ depending on $M$ only such that, for all $\lambda_i \ge 0$ with $\sum_i \lambda_i = 1$,

$$(4.4) \quad \Big|1-P_0\Big(\frac{\sum_i \lambda_i p_{f_i}}{p_{f^*}}\Big)^\alpha - \alpha P_0 \log\frac{p_{f^*}}{\sum_i \lambda_i p_{f_i}}\Big| \le 2\alpha^2 C \sum_i \lambda_i P_0(f_i-f^*)^2.$$

By Lemma 4.3 with $\alpha=1$ and $p=p_f$ and similar arguments, we also have that, for any $f\in\mathscr{F}$,

$$\Big|1-P_0\Big(\frac{\sum_i \lambda_i p_{f_i}}{p_f}\Big) - P_0 \log\frac{p_f}{\sum_i \lambda_i p_{f_i}}\Big| \le 2C \sum_i \lambda_i P_0(f_i-f)^2.$$

For $\lambda_i=1$ this becomes

$$\Big|1-P_0\Big(\frac{p_{f_i}}{p_f}\Big) - P_0 \log\frac{p_f}{p_{f_i}}\Big| \le 2CP_0(f_i-f)^2.$$

Taking differences, we obtain that

$$\Big|P_0 \log\frac{p_f}{\sum_i \lambda_i p_{f_i}} - \sum_i \lambda_i P_0 \log\frac{p_f}{p_{f_i}}\Big| \le 4C\sum_i \lambda_i P_0(f_i-f)^2.$$

By the fact that $\log(ab) = \log a + \log b$ for every $a,b>0$, this inequality remains true if $f$ on the left is replaced by $f^*$. Combine the resulting inequality with (4.4) to find that

$$1-P_0\Big(\frac{\sum_i \lambda_i p_{f_i}}{p_{f^*}}\Big)^\alpha$$
$$\ge \alpha \sum_i \lambda_i P_0 \log\frac{p_{f^*}}{p_{f_i}}$$
$$\quad - 2\alpha^2 C \sum_i \lambda_i P_0(f^*-f_i)^2 - 4C\sum_i \lambda_i P_0(f_i-f)^2$$
$$\ge \Big(\frac{\alpha}{2}-2\alpha^2 C\Big)\sum_i \lambda_i P_0(f^*-f_i)^2 - 4C\sum_i \lambda_i P_0(f_i-f)^2,$$

where we have used (4.2). For sufficiently small $\alpha>0$ and suitable constants $C_2, C_3$, the right-hand side is bounded below by the right-hand side of the lemma. Finally the left-hand side of the lemma can be bounded by the supremum over $\alpha \in (0,1)$ of the left-hand side of the last display, since $-\log x \geq 1 - x$ for every $x > 0$. $\square$

In view of the preceding lemma, the estimation of the quantities involved in the main theorems can be based on the $L_2(P_0)$ distance.

The "neighborhoods" $B(\varepsilon, P_{f^*}; P_0)$ involved in the prior mass conditions (2.4) and (2.11) can be interpreted in the form

$$B(\varepsilon, P_{f^*}; P_0) = \{f \in \mathscr{F} : P_0(f - f_0)^2 \leq P_0(f^* - f_0)^2 + \varepsilon^2, P_0(f - f^*)^2 \leq \varepsilon^2\}.$$

If $P_0(f - f^*)(f^* - f_0) = 0$ for every $f \in \mathscr{F}$ (as is the case if $f^* = f_0$ or if $f^*$ lies in the interior of $\mathscr{F}$), then this reduces to an $L_2(P_0)$-ball around $f^*$ by Pythagoras' theorem.

In view of the preceding lemma and Lemma 2.1, the entropy for testing in (2.5) can be replaced by the local entropy of $\mathscr{F}$ for the $L_2(P_0)$-metric. The rate of convergence of the posterior distribution guaranteed by Theorem 2.1 is then also relative to the $L_2(P_0)$-distance. These observations yield the following theorem.

THEOREM 4.1. *Assume the assertions of Lemma 4.1 and, in addition, that $P_0(f - f^*)(f^* - f_0) = 0$ for every $f \in \mathscr{F}$. If $\varepsilon_n$ is a sequence of strictly positive numbers with $\varepsilon_n \to 0$ and $n\varepsilon_n^2 \to \infty$ such that, for a constant $L > 0$ and all $n$,*

(4.5) $$\Pi(f \in \mathscr{F} : P_0(f - f^*)^2 \leq \varepsilon_n^2) \geq e^{-Ln\varepsilon_n^2},$$

(4.6) $$N(\varepsilon_n, \mathscr{F}, \|\cdot\|_{P_0,2}) \leq e^{n\varepsilon_n^2},$$

*then $\Pi_n(f \in \mathscr{F} : P_0(f - f^*)^2 \geq M\varepsilon_n^2 | X_1, \ldots, X_n) \to 0$ in $L_1(P_0^n)$, for every sufficiently large constant $M$.*

There are many special cases of interest of this theorem and the more general results that can be obtained from Theorems 2.1 and 2.2 using the preceding reasoning. Some of these are considered in the context of the well-specified regression model [14]. The necessary estimates on the prior mass and the entropy are not different for problems other than the regression model. Entropy estimates can also be found in work on rates of convergence of minimum contrast estimators. For these reasons we exclude a discussion of concrete examples.

The following pair of lemmas was used in the proof of the preceding results.

LEMMA 4.2. *There exists a universal constant $C$ such that, for any probability measure $P_0$ and any finite measures $P$ and $Q$ and any $0 < \alpha \leq 1$,*

$$\left|1 - P_0\left(\frac{q}{p}\right)^\alpha - \alpha P_0 \log \frac{p}{q}\right| \leq \alpha^2 C P_0\left[\left(\log \frac{p}{q}\right)^2 \left(\left(\frac{q}{p}\right)^\alpha \mathbb{1}_{\{q > p\}} + \mathbb{1}_{\{q \leq p\}}\right)\right].$$

PROOF. The function $R$ defined by $R(x) = (e^x - 1 - x)/(x^2 e^x)$ for $x \geq 0$ and $R(x) = (e^x - 1 - x)/x^2$ for $x \leq 0$ is uniformly bounded on $\mathbb{R}$ by a constant $C$. We can write

$$P_0 \left(\frac{q}{p}\right)^\alpha = 1 + \alpha P_0 \log \frac{q}{p}$$
$$+ P_0 R\left(\alpha \log \frac{q}{p}\right)\left(\alpha \log \frac{q}{p}\right)^2 \left[\left(\frac{q}{p}\right)^\alpha \mathbb{1}_{\{q>p\}} + \mathbb{1}_{\{q \leq p\}}\right].$$

The lemma follows. □

LEMMA 4.3. *There exists a universal constant $C$ such that, for any probability measure $P_0$ and all finite measures $P, Q_1, \ldots, Q_m$ and constants $0 < \alpha \leq 1$, $\lambda_i \geq 0$ with $\sum_i \lambda_i = 1$,*

$$\left| 1 - P_0 \left(\frac{\sum_i \lambda_i q_i}{p}\right)^\alpha - \alpha P_0 \log \frac{p}{\sum_i \lambda_i q_i} \right| \leq 2\alpha^2 C \sum_i \lambda_i P_0 \left(\log \frac{q_i}{p}\right)^2 \left[\left(\frac{q_i}{p}\right)^2 + 1\right].$$

PROOF. In view of Lemma 4.2 with $q = \sum_i \lambda_i q_i$, it suffices to bound

$$P_0 \left[\left(\log \frac{\sum_i \lambda_i q_i}{p}\right)^2 \left(\left(\frac{\sum_i \lambda_i q_i}{p}\right)^\alpha \mathbb{1}_{\sum_i \lambda_i q_i > p} + \mathbb{1}_{\sum_i \lambda_i q_i \leq p}\right)\right]$$

by the right-hand side of the lemma. We can replace $\alpha$ in the display by 2 and make the expression larger. Next we bound the two terms corresponding to the decomposition by indicators separately.

By the convexity of the map $x \mapsto x \log x$,

$$\left(\log \frac{\sum_i \lambda_i q_i}{p}\right)\left(\frac{\sum_i \lambda_i q_i}{p}\right) \leq \sum_i \lambda_i \left(\log \frac{q_i}{p}\right)\left(\frac{q_i}{p}\right).$$

If $\sum_i \lambda_i q_i > p$, then the left-hand side is positive and the inequality is preserved when we square on both sides. Convexity of the map $x \mapsto x^2$ allows us to bound the square of the right-hand side as in the lemma.

By the concavity of the logarithm,

$$-\log \frac{\sum_i \lambda_i q_i}{p} \leq -\sum_i \lambda_i \log \frac{q_i}{p}.$$

On the set $\sum_i \lambda_i q_i < p$ the left-hand side is positive and we can again take squares on both sides and preserve the inequality. □

4.2. *Laplace regression.* Suppose that the error-density $p$ is equal to the Laplace density $p(x) = \frac{1}{2} \exp(-|x|)$. Then

$$\log \frac{p_f}{p_{f_0}}(X,Y) = -(|e_0 + f_0(X) - f(X)| - |e_0|),$$
$$-P_0 \log \frac{p_f}{p_{f_0}} = P_0 \Phi(f - f_0),$$

for $\Phi(\nu) = \mathrm{E}_0(|e_0 - \nu| - |e_0|)$. The function $\Phi$ is minimized over $\nu \in \mathbb{R}$ at the median of $e_0$. It follows that if $f_0 + m$, for $m$ the median of $e_0$, is contained in $\mathscr{F}$, then the Kullback–Leibler divergence $-P_0 \log(p_f/p_0)$ is minimized over $f \in \mathscr{F}$ at $f = f_0 + m$. If $\mathscr{F}$ is a compact, convex subset of $L_1(P_0)$, then in any case there exists $f^* \in \mathscr{F}$ that minimizes the Kullback–Leibler divergence, but it appears difficult to determine this concretely in general. For simplicity of notation, we shall assume that $m = 0$.

If the distribution of $e_0$ is smooth, then the function $\Phi$ will be smooth too. Because it is minimal at $\nu = m = 0$, it is reasonable to expect that, for $\nu$ in a neighborhood of $m = 0$ and some positive constant $C_0$,

$$(4.7) \qquad \Phi(\nu) = \mathrm{E}_0(|e_0 - \nu| - |e_0|) \geq C_0 |\nu|^2.$$

Because $\Phi$ is convex, it is also reasonable to expect that its second derivative, if it exists, is strictly positive.

LEMMA 4.4. *Let $\mathscr{F}$ be a class of uniformly bounded functions $f : \mathscr{X} \to \mathbb{R}$ and let $f_0$ be uniformly bounded. Assume that either $f_0 \in \mathscr{F}$ and (4.7) holds, or that $\mathscr{F}$ is convex and compact in $L_1(P_0)$ and that $\Phi$ is twice continuously differentiable with strictly positive second derivative. Then there exist positive constants $C_0, C_1, C_2, C_3$ such that, for all $m \in \mathbb{N}$, $f, f_1, \ldots, f_m \in \mathscr{F}$ and $\lambda_1, \ldots, \lambda_m \geq 0$ with $\sum_i \lambda_i = 1$,*

$$P_0 \log \frac{p_f}{p_{f^*}} \leq -C_0 P_0 (f - f^*)^2,$$

$$(4.8) \qquad P_0 \left( \log \frac{p_{f^*}}{p_f} \right)^2 \leq C_1 P_0 (f - f^*)^2,$$

$$\sup_{0 < \alpha < 1} -\log P_0 \left( \frac{\sum_i \lambda_i p_{f_i}}{p_{f^*}} \right)^\alpha \geq C_2 \sum_i \lambda_i (P_0 (f_i - f^*)^2 - C_3 P_0 (f - f_i)^2).$$

PROOF. Suppose first that $f_0 \in \mathscr{F}$, so that $f^* = f_0$. As $\Phi$ is monotone on $(0, \infty)$ and $(-\infty, 0)$, inequality (4.7) is automatically also satisfied for $\nu$ in a given compactum (with $C_0$ depending on the compactum). Choosing the compactum large enough such that $(f - f^*)(X)$ is contained in it with probability one, we conclude that (4.8) holds (with $f_0 = f^*$).

If $f^*$ is not contained in $\mathscr{F}$ but $\mathscr{F}$ is convex, we obtain a similar inequality with $f^*$ replacing $f_0$, as follows. Because $f^*$ minimizes $f \mapsto P_0 \Phi(f - f_0)$ over $\mathscr{F}$ and $f_t = (1-t)f^* + tf \in \mathscr{F}$ for $t \in [0, 1]$, the right derivative of the map $t \mapsto P_0 \Phi(f_t - f_0)$ is nonnegative at $t = 0$. This yields $P_0 \Phi'(f^* - f_0)(f - f^*) \geq 0$. By a Taylor expansion,

$$P_0 \log \frac{p_{f^*}}{p_f} = P_0(\Phi(f - f_0) - \Phi(f^* - f_0))$$

$$= P_0 \Phi'(f^* - f_0)(f - f^*) + \frac{1}{2} P_0 \Phi''(\tilde{f} - f_0)(f - f^*)^2,$$

for some $\tilde{f}$ between $f$ and $f^*$. The first term on the right-hand side is nonnegative and the function $\Phi''$ is bounded away from zero on compacta by assumption. Thus, the right-hand side is bounded below by a constant times $P_0(f - f^*)^2$ and again (4.8) follows.

Because $\log(p_f/p_{f^*})$ is bounded in absolute value by $|f - f^*|$, we also have, with $M$ a uniform upper bound on $\mathscr{F}$ and $f_0$,

$$P_0\left(\log \frac{p_f}{p_{f^*}}\right)^2 \leq P_0(f^* - f)^2,$$

$$P_0\left(\log \frac{p_f}{p_{f^*}}\right)^2 \left(\frac{p_f}{p_{f^*}}\right)^\alpha \leq P_0(f^* - f)^2 e^{2\alpha M}.$$

As in the proof of Lemma 4.1, we can combine these inequalities, (4.8) and Lemma 4.3 to obtain the result. □

As in the case of regression using the normal density for the error-distribution, the preceding lemma reduces the entropy calculations for the application of Theorem 2.1 to estimates of the $L_2(P_0)$-entropy of the class of regression functions $\mathscr{F}$. The resulting rate of convergence is the same as in the case where a normal distribution is used for the error. A difference with the normal case is that presently no tail conditions of the type $E_0 e^{\varepsilon|e_0|} < \infty$ are necessary. Instead the lemma assumes a certain smoothness of the true distribution of the error $e_0$.

**5. Parametric models.** The behavior of posterior distributions for finite-dimensional, misspecified models was considered in [1] and more recently by Bunke and Milhaud [3] (see also the references in the latter). In this section we show that the basic result that the posterior concentrates near a minimal Kullback–Leibler point at the rate $\sqrt{n}$ follows from our general theorems under some natural conditions. We first consider models indexed by a parameter in a general metric space and relate the rate of convergence to the metric entropy of the parameter set. Next we specialize to Euclidean parameter sets.

Let $\{p_\theta : \theta \in \Theta\}$ be a collection of probability densities indexed by a parameter $\theta$ in a metric space $(\Theta, d)$. Let $P_0$ be the true distribution of the data and assume that there exists a $\theta^* \in \Theta$, such that, for all $\theta, \theta_1, \theta_2 \in \Theta$ and some constant $C > 0$,

(5.1) $$P_0 \log \frac{p_\theta}{p_{\theta^*}} \leq -C d^2(\theta, \theta^*),$$

(5.2) $$P_0\left(\sqrt{\frac{p_{\theta_1}}{p_{\theta_2}}} - 1\right)^2 \leq d^2(\theta_1, \theta_2),$$

(5.3) $$P_0\left(\log \frac{p_{\theta_1}}{p_{\theta_2}}\right)^2 \leq d^2(\theta_1, \theta_2).$$

The first inequality implies that $\theta^*$ is a point of minimal Kullback–Leibler divergence $\theta \mapsto -P_0 \log(p_\theta/p_0)$ between $P_0$ and the model. The second and third conditions are (integrated) Lipschitz conditions on the dependence of $p_\theta$ on $\theta$. The following lemma shows that in the application of Theorems 2.1 and 2.2 these conditions allow one to replace the entropy for testing by the local entropy of $\Theta$ relative to (a multiple of) the natural metric $d$.

In examples it may be worthwhile to relax the conditions somewhat. In particular, the conditions (5.2)–(5.3) can be "localized." Rather than assuming that they are valid for every $\theta_1, \theta_2 \in \Theta$, the same results can be obtained if they are valid for every pair $(\theta_1, \theta_2)$ with $d(\theta_1, \theta_2)$ sufficiently small and every pair $(\theta_1, \theta_2)$ with arbitrary $\theta_1$ and $\theta_2 = \theta^*$. For $\theta_2 = \theta^*$ and $P_0 = P_{\theta^*}$ (i.e., the well-specified situation), condition (5.2) is a bound on the Hellinger distance between $P_{\theta^*}$ and $P_{\theta_1}$.

LEMMA 5.1. *Under the preceding conditions, there exist positive constants $C_1, C_2$ such that, for all $m \in \mathbb{N}, \theta, \theta_1, \ldots, \theta_m \in \Theta$ and $\lambda_1, \ldots, \lambda_m \geq 0$ with $\sum_i \lambda_i = 1$,*

$$\sum_i \lambda_i d^2(\theta_i, \theta^*) - C_1 \sum_i \lambda_i d^2(\theta, \theta_i) \leq C_2 \sup_{0 < \alpha < 1} -\log P_0 \left( \frac{\sum_i \lambda_i p_{\theta_i}}{p_{\theta^*}} \right)^\alpha.$$

PROOF. In view of Lemma 5.3 (below) with $p = p_{\theta^*}$, (5.2) and (5.3), there exists a constant $C$ such that

$$(5.4) \quad \left| 1 - P_0 \left( \frac{\sum_i \lambda_i p_{\theta_i}}{p_{\theta^*}} \right)^\alpha - \alpha P_0 \left( \log \frac{p_{\theta^*}}{\sum_i \lambda_i p_{\theta_i}} \right) \right| \leq 2\alpha^2 C \sum_i \lambda_i d^2(\theta_i, \theta^*).$$

By Lemma 5.3 with $\alpha = 1$, $p = p_\theta$, (5.2) and (5.3),

$$\left| 1 - P_0 \left( \frac{\sum_i \lambda_i p_{\theta_i}}{p_\theta} \right) - P_0 \left( \log \frac{p_\theta}{\sum_i \lambda_i p_{\theta_i}} \right) \right| \leq 2C \sum_i \lambda_i d^2(\theta_i, \theta).$$

We can evaluate this with $\lambda_i = 1$ (for each $i$ in turn) and next subtract the convex combination of the resulting inequalities from the preceding display to obtain

$$\left| P_0 \left( \log \frac{p_\theta}{\sum_i \lambda_i p_{\theta_i}} \right) - \sum_i \lambda_i P_0 \left( \log \frac{p_\theta}{p_{\theta_i}} \right) \right| \leq 4C \sum_i \lambda_i d^2(\theta_i, \theta).$$

By the additivity of the logarithm, this remains valid if we replace $\theta$ on the left-hand side by $\theta^*$. Combining the resulting inequality with (5.1) and (5.4), we obtain

$$1 - P_0 \left( \frac{\sum_i \lambda_i p_{\theta_i}}{p_{\theta^*}} \right)^\alpha \geq \alpha \sum_i \lambda_i d^2(\theta_i, \theta^*)(C - 2\alpha) - 4C \sum_i \lambda_i d^2(\theta_i, \theta).$$

The lemma follows upon choosing $\alpha > 0$ sufficiently small and using $-\log x \geq 1 - x$. □

If the prior on the model $\{p_\theta : \theta \in \Theta\}$ is induced by a prior on the parameter set $\Theta$, then the prior mass condition (2.11) translates into a lower bound for the prior mass of the set

$$B(\varepsilon, \theta^*; P_0) = \left\{\theta \in \Theta : -P_0 \log \frac{p_\theta}{p_{\theta^*}} \leq \varepsilon^2, P_0\left(\log \frac{p_\theta}{p_{\theta^*}}\right)^2 \leq \varepsilon^2\right\}.$$

In addition to (5.1), it is reasonable to assume a lower bound of the form

$$(5.5) \qquad P_0 \log \frac{p_\theta}{p_{\theta^*}} \geq -\underline{C} d^2(\theta, \theta^*),$$

at least for small values of $d(\theta, \theta^*)$. This together with (5.3) implies that $B(\varepsilon, \theta^*; P_0)$ contains a ball of the form $\{\theta : d(\theta, \theta^*) \leq C_1 \varepsilon\}$ for small enough $\varepsilon$. Thus, in the verification of (2.4) or (2.11) we may replace $B(\varepsilon, P^*; P_0)$ by a ball of radius $\varepsilon$ around $\theta^*$. These observations lead to the following theorem.

THEOREM 5.1. *Let* (5.1)–(5.5) *hold. If for sufficiently small $A$ and $C$,*

$$\sup_{\varepsilon > \varepsilon_n} \log N(A\varepsilon, \{\theta \in \Theta : \varepsilon < d(\theta, \theta^*) < 2\varepsilon\}, d) \leq n\varepsilon_n^2,$$

$$\frac{\Pi(\theta : j\varepsilon_n < d(\theta, \theta^*) < 2j\varepsilon_n)}{\Pi(\theta : d(\theta, \theta^*) \leq C\varepsilon_n)} \leq e^{n\varepsilon_n^2 j^2/8},$$

*then* $\Pi(\theta : d(\theta, \theta^*) \geq M_n \varepsilon_n | X_1, \ldots, X_n) \to 0$ *in* $L_1(P_0^n)$ *for any* $M_n \to \infty$.

5.1. *Finite-dimensional models.* Let $\Theta$ be an open subset of $m$-dimensional Euclidean space equipped with the Euclidean distance $d$ and let $\{p_\theta : \theta \in \Theta\}$ be a model satisfying (5.1)–(5.5).

Then the local covering numbers as in the preceding theorem satisfy, for some constant $B$,

$$N(A\varepsilon, \{\theta \in \Theta : \varepsilon < d(\theta, \theta^*) < 2\varepsilon\}, d) \leq \left(\frac{B}{A}\right)^m$$

(see, e.g., [8], Section 5). In view of Lemma 2.2, condition (2.5) is satisfied for $\varepsilon_n$ a large multiple of $1/\sqrt{n}$. If the prior $\Pi$ on $\Theta$ possesses a density that is bounded away from zero and infinity, then

$$\frac{\Pi(\theta : d(\theta, \theta^*) \leq j\varepsilon)}{\Pi(B(\varepsilon, \theta^*; P_0))} \leq C_2 j^m,$$

for some constant $C_2$. It follows that (2.11) is satisfied for the same $\varepsilon_n$. Hence, the posterior concentrates at rate $1/\sqrt{n}$ near the point $\theta^*$.

The preceding situation arises if the minimal point $\theta^*$ is interior to the parameter set $\Theta$. An example is fitting an exponential family, such as the Gaussian model, to observations that are not sampled from an element of the family. If the minimal point $\theta^*$ is not interior to $\Theta$, then we cannot expect (5.1) to hold for the natural distance and different rates of convergence may arise. We include a simple example of the latter type, which is somewhat surprising.

EXAMPLE. Suppose that $P_0$ is the standard normal distribution and the model consists of all $N(\theta, 1)$-distributions with $\theta \geq 1$. The minimal Kullback–Leibler point is $\theta^* = 1$. If the prior possesses a density on $[1, \infty)$ that is bounded away from 0 and infinity near 1, then the posterior concentrates near $\theta^*$ at the rate $1/n$.

One easily shows that

$$
\begin{aligned}
-P_0 \log \frac{p_\theta}{p_{\theta^*}} &= \frac{1}{2}(\theta - \theta^*)(\theta + \theta^*), \\
-\log P_0 \left(\frac{p_\theta}{p_{\theta^*}}\right)^\alpha &= \frac{1}{2}\alpha(\theta - \theta^*)(\theta + \theta^* - \alpha(\theta - \theta^*)).
\end{aligned}
\tag{5.6}
$$

This shows that (2.9) is satisfied for a multiple of the metric $d(p_{\theta_1}, p_{\theta_2}) = \sqrt{|\theta_1 - \theta_2|}$ on $\Theta = [1, \infty)$. Its strengthening (2.8) can be verified by the same methods as before, or, alternatively, the existence of suitable tests can be established directly based on the special nature of the normal location family. [A suitable test for an interval $(\theta_1, \theta_2)$ can be obtained from a suitable test for its left end point.] The entropy and prior mass can be estimated as in regular parametric models and conditions (2.5)–(2.11) can be shown to be satisfied for $\varepsilon_n$ a large multiple of $1/\sqrt{n}$. This yields the rate $1/\sqrt{n}$ relative to the metric $\sqrt{|\theta_1 - \theta_2|}$ and, hence, the rate $1/n$ in the natural metric.

Theorem 2.2 only gives an upper bound on the rate of convergence. In the present situation this appears to be sharp. For instance, for a uniform prior on $[1, 2]$, the posterior mass of the interval $[c, 2]$ can be seen to be, with $Z_n = \sqrt{n} \bar{X}_n$,

$$
\frac{\Phi(2\sqrt{n} - Z_n) - \Phi(c\sqrt{n} - Z_n)}{\Phi(2\sqrt{n} - Z_n) - \Phi(\sqrt{n} - Z_n)} \approx \frac{\sqrt{n} - Z_n}{c\sqrt{n} - Z_n} e^{(-1/2)(c^2 - 1)n + Z_n(c-1)\sqrt{n}},
$$

where we use Mills' ratio to see that $\Phi(y_n) - \Phi(x_n) \approx (1/x_n)\phi(x_n)$ if $x_n, y_n \to c \in (0, 1)$ such that $x_n/y_n \to 0$. This is bounded away from zero for $c = c_n = 1 + C/n$ and fixed $C$.

LEMMA 5.2. *There exists a universal constant $C$ such that for any probability measure $P_0$ and any finite measures $P$ and $Q$ and any $0 < \alpha \leq 1$,*

$$
\left|1 - P_0\left(\frac{q}{p}\right)^\alpha - \alpha P_0 \log \frac{p}{q}\right| \leq \alpha^2 C P_0\left[\left(\sqrt{\frac{q}{p}} - 1\right)^2 \mathbb{1}_{\{q > p\}} + \left(\log \frac{p}{q}\right)^2 \mathbb{1}_{\{q \leq p\}}\right].
$$

LEMMA 5.3. *There exists a universal constant $C$ such that, for any probability measure $P_0$ and any finite measures $P$, $Q_1, \ldots, Q_m$ and any $\lambda_1, \ldots, \lambda_m \geq 0$ with $\sum_i \lambda_i = 1$ and $0 < \alpha \leq 1$, the following inequality holds:*

$$\left| 1 - P_0 \left( \frac{\sum_i \lambda_i q_i}{p} \right)^\alpha - \alpha P_0 \log \frac{p}{\sum_i \lambda_i q_i} \right|$$
$$\leq 2\alpha^2 C \sum_i \lambda_i P_0 \left[ \left( \sqrt{\frac{q_i}{p}} - 1 \right)^2 + \left( \log \frac{q_i}{p} \right)^2 \right].$$

PROOFS. The function $R$ defined by $R(x) = (e^x - 1 - x)/\alpha^2(e^{x/2\alpha} - 1)^2$ for $x \geq 0$ and $R(x) = (e^x - 1 - x)/x^2$ for $x \leq 0$ is uniformly bounded on $\mathbb{R}$ by a constant $C$, independent of $\alpha \in (0, 1]$. [This may be proved by noting that the functions $(e^x - 1)/\alpha(e^{\alpha x} - 1)$ and $(e^x - 1 - x)/(e^{x/2} - 1)^2$ are bounded, where this follows for the first by developing the exponentials in their power series.] For the proof of the first lemma, we can proceed as in the proof of Lemma 4.2. For the proof of the second lemma, we proceed as in the proof of Lemma 4.3, this time also making use of the convexity of the map $x \mapsto |\sqrt{x} - 1|^2$ on $[0, \infty)$. □

**6. Existence of tests.** The proofs of Theorems 2.1 and 2.2 rely on tests of $P_0$ versus the positive, finite measures $Q(P)$ obtained from points $P$ that are at positive distance from the set of points with minimal Kullback–Leibler divergence. Because we need to test $P_0$ against finite measures (i.e., not necessarily probability measures), known results on tests using the Hellinger distance, such as in [10] or [8], do not apply. It turns out that in this situation the Hellinger distance may not be appropriate and instead we use the full Hellinger transform. The aim of this section is to prove the existence of suitable tests and give upper bounds on their power. We first formulate the results in a general notation and then specialize to the application in misspecified models.

6.1. *General setup.* Let $P$ be a probability measure on a measurable space $(\mathscr{X}, \mathscr{A})$ (playing the role of $P_0$) and let $\mathscr{Q}$ be a class of finite measures on $(\mathscr{X}, \mathscr{A})$ [playing the role of the measures $Q$ with $dQ = (p_0/p^*)\, dP$]. We wish to bound the minimax risk for testing $P$ versus $\mathscr{Q}$, defined by

$$\pi(P, \mathscr{Q}) = \inf_\phi \sup_{Q \in \mathscr{Q}} (P\phi + Q(1 - \phi)),$$

where the infimum is taken over all measurable functions $\phi : \mathscr{X} \to [0, 1]$. Let $\operatorname{conv}(\mathscr{Q})$ denote the convex hull of the set $\mathscr{Q}$.

LEMMA 6.1. *If there exists a $\sigma$-finite measure that dominates all $Q \in \mathcal{Q}$, then*

$$\pi(P, \mathcal{Q}) = \sup_{Q \in \mathrm{conv}(\mathcal{Q})} (P(p < q) + Q(p \geq q)).$$

*Moreover, there exists a test $\phi$ that attains the infimum in the definition of $\pi(P, \mathcal{Q})$.*

PROOF. If $\mu'$ is a measure dominating $\mathcal{Q}$, then a $\sigma$-finite measure $\mu$ exists that dominates both $\mathcal{Q}$ and $P$ (e.g., $\mu = \mu' + P$). Let $p$ and $q$ be $\mu$-densities of $P$ and $Q$, for every $Q \in \mathcal{Q}$. The set of test-functions $\phi$ can be identified with the positive unit ball $\Phi$ of $L_\infty(\mathcal{X}, \mathcal{A}, \mu)$, which is dual to $L_1(\mathcal{X}, \mathcal{A}, \mu)$, since $\mu$ is $\sigma$-finite. If equipped with the weak-$*$ topology, the positive unit ball $\Phi$ is Hausdorff and compact by the Banach–Alaoglu theorem (see, e.g., [11], Theorem 2.6.18, and note that the positive functions form a closed and convex subset of the unit ball). The convex hull $\mathrm{conv}(\mathcal{Q})$ (or rather the corresponding set of $\mu$-densities) is a convex subset of $L_1(\mathcal{X}, \mathcal{A}, \mu)$. The map

$$L_\infty(\mathcal{X}, \mathcal{A}, \mu) \times L_1(\mathcal{X}, \mathcal{A}, \mu) \to \mathbb{R},$$
$$(\phi, Q) \mapsto \phi P + (1 - \phi) Q$$

is concave in $Q$ and convex in $\phi$. [Note that in the current context we write $\phi P$ instead of $P\phi$, in accordance with the fact that we consider $\phi$ as a bounded linear functional on $L_1(\mathcal{X}, \mathcal{A}, \mu)$.] Moreover, the map is weak-$*$-continuous in $\phi$ for every fixed $Q$ [note that every weak-$*$-converging net $\phi_\alpha \xrightarrow{\text{w-}*} \phi$ by definition satisfies $\phi_\alpha Q \to \phi Q$ for all $Q \in L_1(\mathcal{X}, \mathcal{A}, \mu)$]. The conditions for application of the minimax theorem (see, e.g., [15], page 239) are satisfied and we conclude

$$\inf_{\phi \in \Phi} \sup_{Q \in \mathrm{conv}(\mathcal{Q})} (\phi P + (1 - \phi) Q) = \sup_{Q \in \mathrm{conv}(\mathcal{Q})} \inf_{\phi \in \Phi} (\phi P + (1 - \phi) Q).$$

The expression on the left-hand side is the minimax testing risk $\pi(P, \mathcal{Q})$. The infimum on the right-hand side is attained at the point $\phi = \mathbb{1}\{p < q\}$, which leads to the first assertion of the lemma upon substitution.

The second assertion of the lemma follows because the function $\phi \mapsto \sup\{\phi P + (1 - \phi) Q : Q \in \mathrm{conv}(\mathcal{Q})\}$ is a supremum of weak-$*$-continuous functions and, hence, attains its minimum on the compactum $\Phi$. □

It is possible to express the right-hand side of the preceding lemma in the $L_1$-distance between $P$ and $Q$, but this is not useful for the following. Instead, we use a bound in terms of the *Hellinger transform* $\rho_\alpha(P, Q)$ defined by, for $0 < \alpha < 1$,

$$\rho_\alpha(P, Q) = \int p^\alpha q^{1-\alpha} \, d\mu.$$

By Hölder's inequality, this quantity is finite for all finite measures $P$ and $Q$. The definition is independent of the choice of dominating measure $\mu$.

For any pair $(P,Q)$ and every $\alpha \in (0,1)$, we can bound

$$
\begin{aligned}
P(p<q) + Q(p\geq q) &= \int_{p<q} p\,d\mu + \int_{p\geq q} q\,d\mu \\
&\leq \int_{p<q} p^\alpha q^{1-\alpha}\,d\mu + \int_{p\geq q} p^\alpha q^{1-\alpha}\,d\mu \\
&= \rho_\alpha(P,Q).
\end{aligned}
\tag{6.1}
$$

Hence, the right-hand side of the preceding lemma is bounded by $\sup_Q \rho_\alpha(P,Q)$ for all $\alpha \in (0,1)$. The advantage of this bound is the fact that it factorizes if $P$ and $Q$ are product measures. For ease of notation, define

$$\rho_\alpha(\mathscr{P},\mathscr{Q}) = \sup\{\rho_\alpha(P,Q): P \in \mathrm{conv}(\mathscr{P}), Q \in \mathrm{conv}(\mathscr{Q})\}.$$

LEMMA 6.2. *For any $0 < \alpha < 1$ and classes $\mathscr{P}_1, \mathscr{P}_2, \mathscr{Q}_1, \mathscr{Q}_2$ of finite measures,*

$$\rho_\alpha(\mathscr{P}_1 \times \mathscr{P}_2, \mathscr{Q}_1 \times \mathscr{Q}_2) \leq \rho_\alpha(\mathscr{P}_1, \mathscr{Q}_1)\rho_\alpha(\mathscr{P}_2, \mathscr{Q}_2),$$

*where $\mathscr{P}_1 \times \mathscr{P}_2$ denotes the class of product measures $\{P_1 \times P_2 : P_1 \in \mathscr{P}_1, P_2 \in \mathscr{P}_2\}$.*

PROOF. Let $P \in \mathrm{conv}(\mathscr{P}_1 \times \mathscr{P}_2)$ and $Q \in \mathrm{conv}(\mathscr{Q}_1 \times \mathscr{Q}_2)$ be given. Since both are (finite) convex combinations, $\sigma$-finite measures $\mu_1$ and $\mu_2$ can always be found such that both $P$ and $Q$ have $\mu_1 \times \mu_2$ densities which both can be written in the form of a finite convex combination as

$$p(x,y) = \sum_i \lambda_i p_{1i}(x) p_{2i}(y), \qquad \lambda_i \geq 0, \ \sum_i \lambda_i = 1,$$

$$q(x,y) = \sum_j \kappa_j q_{1j}(x) q_{2j}(y), \qquad \kappa_j \geq 0, \ \sum_j \kappa_j = 1,$$

for $\mu_1 \times \mu_2$-almost-all pairs $(x,y) \in \mathscr{X} \times \mathscr{X}$. Here $p_{1i}$ and $q_{1j}$ are $\mu_1$-densities for measures belonging to $\mathscr{P}_1$ and $\mathscr{Q}_1$, respectively (and, analogously, $p_{2i}$ and $q_{2j}$ are $\mu_2$-densities for measures in $\mathscr{P}_2$ and $\mathscr{Q}_2$). This implies that we can write

$$\int p^\alpha q^{1-\alpha}\,d(\mu_1 \times \mu_2)$$

$$= \int \left\{ \int \left(\frac{\sum_i \lambda_i p_{1i}(x) p_{2i}(y)}{\sum_i \lambda_i p_{1i}(x)}\right)^\alpha \left(\frac{\sum_j \kappa_j q_{1j}(x) q_{2j}(y)}{\sum_j \kappa_j q_{1j}(x)}\right)^{1-\alpha} d\mu_2(y) \right\}$$

$$\times \left(\sum_i \lambda_i p_{1i}(x)\right)^\alpha \left(\sum_j \kappa_j q_{1j}(x)\right)^{1-\alpha} d\mu_1(x)$$

(where, as usual, the integrand of the inner integral is taken equal to zero whenever the $\mu_1$-density equals zero). The inner integral is bounded by $\rho_\alpha(\mathscr{P}_2, \mathscr{Q}_2)$ for every fixed $x \in \mathscr{X}$. After substituting this upper bound, the remaining integral is bounded by $\rho_\alpha(\mathscr{P}_1, \mathscr{Q}_1)$. $\square$

Combining (6.1) with Lemmas 6.2 and 6.1, we obtain the following theorem.

THEOREM 6.1. *If $P$ is a probability measure on $(\mathscr{X}, \mathscr{A})$ and $\mathscr{Q}$ is a dominated set of finite measures on $(\mathscr{X}, \mathscr{A})$, then, for every $n \geq 1$, there exists a test $\phi_n : \mathscr{X}^n \to [0,1]$ such that, for all $0 < \alpha < 1$,*

$$\sup_{Q \in \mathscr{Q}} (P^n \phi_n + Q^n(1 - \phi_n)) \leq \rho_\alpha(P, \mathscr{Q})^n.$$

The bound given by the theorem is useful only if $\rho_\alpha(P, \mathscr{Q}) < 1$. For probability measures $P$ and $Q$, we have

$$\rho_{1/2}(P, Q) = 1 - \tfrac{1}{2} \int (\sqrt{p} - \sqrt{q})^2 \, d\mu$$

and, hence, we might use the bound with $\alpha = 1/2$ if the Hellinger distance of $\operatorname{conv}(\mathscr{Q})$ to $P$ is positive. For a general finite measure $Q$, the quantity $\rho_{1/2}(P, Q)$ may be bigger than 1 and, depending on $Q$, the Hellinger transform $\rho_\alpha(P, Q)$ may even lie above 1 for every $\alpha$. The following lemma shows that this is controlled by the (generalized) Kullback–Leibler divergence $-P \log(q/p)$.

LEMMA 6.3. *For a probability measure $P$ and a finite measure $Q$, the function $\alpha \mapsto \rho_\alpha(Q, P)$ is convex on $[0,1]$ with $\rho_\alpha(Q, P) \to P(q > 0)$ as $\alpha \downarrow 0$, $\rho_\alpha(Q, P) \to Q(p > 0)$ as $\alpha \uparrow 1$ and*

$$\left. \frac{d\rho_\alpha(Q, P)}{d\alpha} \right|_{\alpha=0} = P \log\left(\frac{q}{p}\right) \mathbb{1}_{q>0}$$

*(which may be equal to $-\infty$).*

PROOF. The function $\alpha \mapsto e^{\alpha y}$ is convex on $(0,1)$ for all $y \in [-\infty, \infty)$, implying the convexity of $\alpha \mapsto \rho_\alpha(Q, P) = P(q/p)^\alpha$ on $(0,1)$. The function $\alpha \mapsto y^\alpha = e^{\alpha \log y}$ is continuous on $[0,1]$ for any $y > 0$, is decreasing for $y < 1$, increasing for $y > 1$ and constant for $y = 1$. By monotone convergence, as $\alpha \downarrow 0$,

$$Q\left(\frac{p}{q}\right)^\alpha \mathbb{1}_{\{0<p<q\}} \uparrow Q\left(\frac{p}{q}\right)^0 \mathbb{1}_{\{0<p<q\}} = Q(0 < p < q).$$

By the dominated convergence theorem, with dominating function $(p/q)^\alpha \times \mathbb{1}_{\{p \geq q\}} \leq (p/q)^{1/2} \mathbb{1}_{\{p \geq q\}}$ for $\alpha \leq 1/2$, we have, as $\alpha \to 0$,

$$Q\left(\frac{p}{q}\right)^\alpha \mathbb{1}_{\{p \geq q\}} \to Q\left(\frac{p}{q}\right)^0 \mathbb{1}_{\{p \geq q\}} = Q(p \geq q).$$

Combining the two preceding displays above, we see that $\rho_{1-\alpha}(Q, P) = Q(p/q)^\alpha \to Q(p > 0)$ as $\alpha \downarrow 0$.

By the convexity of the function $\alpha \mapsto e^{\alpha y}$, the map $\alpha \mapsto f_\alpha(y) = (e^{\alpha y} - 1)/\alpha$ decreases, as $\alpha \downarrow 0$, to $(d/d\alpha)|_{\alpha=0} f_\alpha(y) = y$, for every $y$. For $y \leq 0$, we have $f_\alpha(y) \leq 0$, while, for $y \geq 0$, by Taylor's formula,

$$f_\alpha(y) \leq \sup_{0 < \alpha' \leq \alpha} y e^{\alpha' y} \leq y e^{\alpha y} \leq \frac{1}{\varepsilon} e^{(\alpha + \varepsilon) y}.$$

Hence, we conclude that $f_\alpha(y) \leq 0 \vee \varepsilon^{-1} e^{(\alpha+\varepsilon)y} \mathbb{1}_{y \geq 0}$. Consequently, the quotient $\alpha^{-1}(e^{\alpha \log(q/p)} - 1)$ decreases to $\log(q/p)$ as $\alpha \downarrow 0$ and is bounded above by $0 \vee \varepsilon^{-1}(q/p)^{2\varepsilon} \mathbb{1}_{q \geq p}$ for small $\alpha > 0$, which is $P$-integrable for $2\varepsilon < 1$. We conclude that

$$\frac{1}{\alpha}(\rho_\alpha(Q, P) - \rho_0(Q, P)) = \frac{1}{\alpha} P\left(\left(\frac{q}{p}\right)^\alpha - 1\right) \mathbb{1}_{q > 0} \downarrow P \log\left(\frac{q}{p}\right) \mathbb{1}_{q > 0},$$

as $\alpha \downarrow 0$, by the monotone convergence theorem. $\square$

Two typical graphs of the Hellinger transform $\alpha \mapsto \rho_\alpha(Q, P)$ are shown in Figure 1 [corresponding to fitting a unit variance normal location model in a situation where the observations are sampled from a $N(0, 2)$-distribution]. For $P$ a probability measure with $P \ll Q$, the Hellinger transform is equal to 1 at $\alpha = 0$, but will eventually increase to a level that is above 1 near $\alpha = 1$ if $Q(p > 0) > 1$. Unless the slope $P \log(p/q)$ is negative, it will never decrease below the level 1. For probability measures $P$ and $Q$, this slope equals minus the Kullback–Leibler distance and, hence, is strictly negative unless $P = Q$. In that case, the graph is strictly below 1 on $(0, 1)$ and $\rho_{1/2}(P, Q)$ is a convenient choice to work with. For a general finite measure $Q$, the Hellinger transform $\rho_\alpha(Q, P)$ is guaranteed to assume values strictly less than 1 near $\alpha = 0$, provided that the Kullback–Leibler divergence $P \log(p/q)$ is negative, which is not automatically the case. For testing a composite alternative $\mathscr{Q}$, we shall need that this is the case uniformly in $Q \in \text{conv}(\mathscr{Q})$. For a convex alternative $\mathscr{Q}$, Theorem 6.1 guarantees the existence of tests based on $n$ observations with error probabilities bounded by $e^{-n\varepsilon^2}$ if

$$\varepsilon^2 \leq \sup_{0 < \alpha < 1} \sup_{Q \in \mathscr{Q}} \log \frac{1}{\rho_\alpha(Q, P)}.$$

In some of the examples we can achieve inequalities of this type by bounding the right-hand side below by a (uniform) Taylor expansion of $\alpha \mapsto$ $-\log \rho_\alpha(P,Q)$ in $\alpha$ near $\alpha = 0$. Such arguments are not mere technical generalizations: they can be necessary already to prove posterior consistency relative to misspecified standard parametric models.

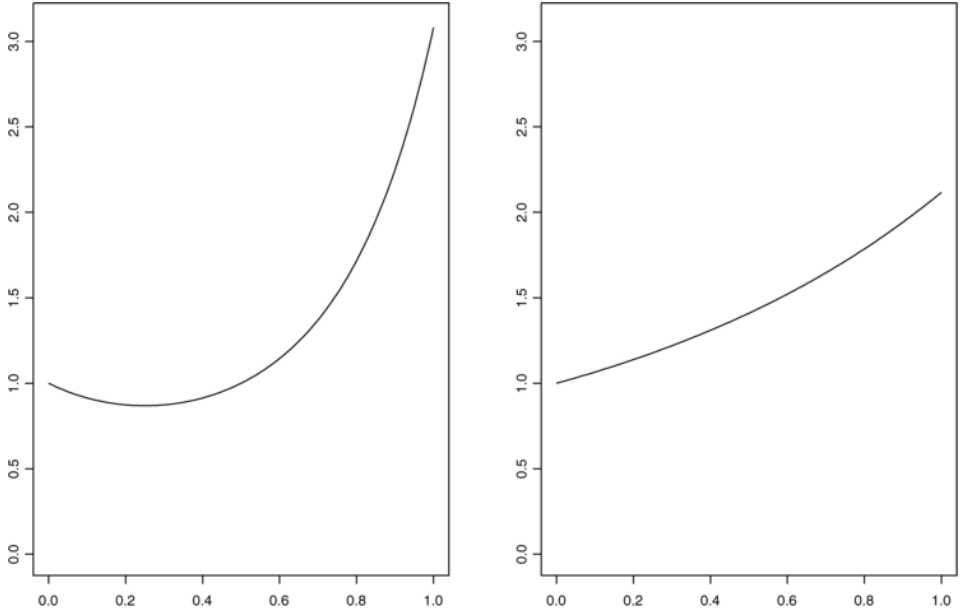

Fig. 1. *The Hellinger transforms* $\alpha \mapsto \rho_\alpha(Q,P)$, *for* $P = N(0,2)$ *and* $Q$, *respectively, the measure defined by* $dQ = (dN(3/2,1)/dN(0,1))\,dP$ (left) *and* $dQ = (dN(3/2,1)/dN(1,1))\,dP$ (right). *Intercepts with the vertical axis at the right and left of the graphs equal* $P(q>0)$ *and* $Q(p>0)$, *respectively. The slope at* 0 *equals (minus) the Kullback–Leibler divergence* $P\log(p/q)$.

If $P(q=0) > 0$, then the Hellinger transform is strictly less than 1 at $\alpha = 0$ and, hence, good tests exist, even though it may be true that $\rho_{1/2}(P,Q) > 1$. The existence of good tests is obvious in this case, since we can reject $Q$ if the observations land in the set $q = 0$.

In the above we have assumed that $\mathscr{Q}$ is dominated. If this is not the case, then the results go through, provided that we use Le Cam's generalized tests [10], that is, we define

$$\pi(P,\mathscr{Q}) = \inf_\phi \sup_{Q \in \mathscr{Q}} (\phi P + (1-\phi)Q),$$

where the infimum is taken over the set of all continuous, positive linear maps $\phi : L_1(\mathscr{X},\mathscr{A}) \mapsto \mathbb{R}$ such that $\phi P \le 1$ for all probability measures $P$. This collection of functionals includes the linear maps that arise from integration of measurable functions $\phi : \mathscr{X} \mapsto [0,1]$, but may be larger. Such tests would be good enough for our purposes, but the generality appears to have little additional value for our application to misspecified models.

The next step is to extend the upper bound to alternatives $\mathscr{Q}$ that are possibly not convex. We are particularly interested in alternatives that are complements of balls around $P$ in some metric. Let $L_1^+(\mathscr{X}, \mathscr{A})$ be the set of finite measures on $(\mathscr{X}, \mathscr{A})$ and let $\tau: L_1^+(\mathscr{X}, \mathscr{A}) \times L_1^+(\mathscr{X}, \mathscr{A}) \mapsto \mathbb{R}$ be such that $\tau(P, \cdot): \mathscr{Q} \mapsto \mathbb{R}$ is a nonnegative function (written in a notation so as to suggest a distance from $P$ to $Q$). For $Q \in \mathscr{Q}$, set

$$\bar{\tau}^2(P, Q) = \sup_{0 < \alpha < 1} \log \frac{1}{\rho_\alpha(P, Q)}. \tag{6.2}$$

For $\varepsilon > 0$, define $N_\tau(\varepsilon, \mathscr{Q})$ to be the minimal number of convex subsets of $\{Q \in L_1^+(\mathscr{X}, \mathscr{A}): \bar{\tau}(P, Q) > \varepsilon/2\}$ needed to cover $\{Q \in \mathscr{Q}: \varepsilon < \tau(P, Q) < 2\varepsilon\}$ and assume that $\mathscr{Q}$ is such that this number is finite for all $\varepsilon > 0$. (The requirement that these convex subsets have $\bar{\tau}$-distance $\varepsilon/2$ to $P$ is essential.) Then the following theorem applies.

THEOREM 6.2. *Let $P$ be a probability measure and $\mathscr{Q}$ be a dominated set of finite measures on $(\mathscr{X}, \mathscr{A})$. Then for all $\varepsilon > 0$ and all $n \geq 1$, there exists a test $\phi_n$ such that, for all $J \in \mathbb{N}$,*

$$P^n \phi_n \leq \sum_{j=1}^\infty N_\tau(j\varepsilon, \mathscr{Q}) e^{-nj^2 \varepsilon^2/4},$$

(6.3)

$$\sup_{\{Q: \tau(P,Q) > J\varepsilon\}} Q^n(1 - \phi_n) \leq e^{-nJ^2 \varepsilon^2/4}.$$

PROOF. Fix $n \geq 1$ and $\varepsilon > 0$ and define $\mathscr{Q}_j = \{Q \in \mathscr{Q}: j\varepsilon < \tau(P, Q) \leq (j+1)\varepsilon\}$. By assumption, there exists for every $j \geq 1$ a finite cover of $\mathscr{Q}_j$ by $N_j = N_\tau(j\varepsilon, \mathscr{Q})$ convex sets $C_{j,1}, \ldots, C_{j,N_j}$ of finite measures, with the further property that

$$\inf_{Q \in C_{j,i}} \bar{\tau}(P, Q) > \frac{j\varepsilon}{2}, \qquad 1 \leq i \leq N_j. \tag{6.4}$$

According to Theorem 6.1, for all $n \geq 1$ and for each set $C_{j,i}$, there exists a test $\phi_{n,j,i}$ such that, for all $\alpha \in (0, 1)$, we have

$$P^n \phi_{n,j,i} \leq \rho_\alpha(P, C_{j,i})^n,$$

$$\sup_{Q \in C_{j,i}} Q^n(1 - \phi_{n,j,i}) \leq \rho_\alpha(P, C_{j,i})^n.$$

By (6.4), we have

$$\sup_{Q \in C_{j,i}} \inf_{0 < \alpha < 1} \rho_\alpha(P, Q) = \sup_{Q \in C_{j,i}} e^{-\bar{\tau}^2(P, Q)} \leq e^{-j^2 \varepsilon^2/4}.$$

For fixed $P$ and $Q$, the function $\alpha \mapsto \rho_\alpha(P, Q)$ is convex and can be extended continuously to a convex function on $[0, 1]$. The function $Q \mapsto \rho_\alpha(P, Q)$ with domain $L_1^+(\mathscr{X}, \mathscr{A})$ is concave. By the minimax theorem (see, e.g., [15], page 239), the left-hand side of the preceding display equals

$$\inf_{0<\alpha<1} \sup_{Q \in C_{j,i}} \rho_\alpha(P, Q) = \inf_{0<\alpha<1} \rho_\alpha(P, C_{j,i}).$$

It follows that

$$P^n \phi_{n,j,i} \vee \sup_{Q \in C_{j,i}} Q^n(1 - \phi_{n,j,i}) \leq e^{-nj^2 \varepsilon^2/4}.$$

Now define a new test function $\phi_n$ by

$$\phi_n = \sup_{j \geq 1} \max_{1 \leq i \leq N_j} \phi_{n,j,i}.$$

Then, for every $J \geq 1$,

$$P^n \phi_n \leq \sum_{j=1}^\infty \sum_{i=1}^{N_j} P^n \phi_{n,j,i} \leq \sum_{j=1}^\infty N_j e^{-nj^2 \varepsilon^2/4},$$

$$\sup_{Q \in \mathscr{Q}} Q^n(1 - \phi_n) \leq \sup_{j \geq J} \max_{i \leq N_j} \sup_{Q \in C_{j,i}} Q^n(1 - \phi_{n,j,i}) \leq \sup_{j \geq J} e^{-nj^2 \varepsilon^2/4} = e^{-nJ^2 \varepsilon^2/4},$$

where $\mathscr{Q} = \{Q : \tau(P, Q) > J\varepsilon\} = \bigcup_{j \geq J} \mathscr{Q}_j$. □

6.2. *Application to misspecification.* When applying the above in the proof for consistency in misspecified models, the problem is to test the true distribution $P_0$ against measures $Q = Q(P)$ taking the form $dQ = (p_0/p^*) \, dP$ for $P \in \mathscr{P}$. In this case, the Hellinger transform takes the form $\rho_\alpha(Q, P_0) = P_0(p/p^*)^\alpha$ and its right derivative at $\alpha = 0$ is equal to $P_0 \log(p/p^*)$. This is negative for every $P \in \mathscr{P}$ if and only if $P^*$ is the point in $\mathscr{P}$ at minimal Kullback–Leibler divergence to $P_0$. This observation illustrates that the measure $P^*$ in Theorem 2.1 is necessarily a point of minimal Kullback–Leibler divergence, even if this is not explicitly assumed. We formalize this in the following lemma.

LEMMA 6.4. *If $P^*$ and $P$ are such that $P_0 \log(p_0/p^*) < \infty$ and $P_0(p/p^*) < \infty$ and the right-hand side of (2.9) is positive, then $P_0 \log(p_0/p^*) < P_0 \log(p_0/p)$. Consequently, the covering numbers for testing $N_t(\varepsilon, \mathscr{P}, d; P_0, P^*)$ in Theorem 2.1 can be finite only if $P^*$ is a point of minimal Kullback–Leibler divergence relative to $P_0$.*

PROOF. The assumptions imply that $P_0(p^* > 0) = 1$. If $P_0(p = 0) > 0$, then $P_0 \log(p_0/p) = \infty$ and there is nothing to prove. Thus, we may assume that $p$ is also strictly positive under $P_0$. Then, in view of Lemma 6.3, the function $g$ defined by $g(\alpha) = P_0(p/p^*)^\alpha = \rho_\alpha(Q, P_0)$ is continuous on $[0, 1]$ with $g(0) = P_0(p > 0) = 1$ and the right-hand side of (2.9) can be positive only if $g(\alpha) < 1$ for some $\alpha \in [0, 1]$. By convexity of $g$ and the fact that $g(0) = 1$, this can happen only if the right derivative of $g$ at zero is nonpositive. In view of Lemma 6.3, this derivative is $g'(0+) = P_0 \log(p/p^*)$.

Finiteness of the covering numbers for testing for some $\varepsilon > 0$ implies that the right-hand side of (2.9) is positive every $P \in \mathscr{P}$ with $d(P, P^*) > 0$, as every such $P$ must be contained in one of the sets $B_i$ in the definition of $N_t(\varepsilon, \mathscr{P}, d; P_0, P^*)$ for some $\varepsilon > 0$, in which case the right-hand side of (2.9) is bounded below by $\varepsilon^2/4$.

If $P_0(p/p^*) \leq 1$ for every $P \in \mathscr{P}$, then the measure $Q$ defined by $dQ = (p_0/p^*) \, dP$ is a subprobability measure and, hence, by convexity, the Hellinger transform $\alpha \mapsto \rho_\alpha(P_0, Q)$ is never above the level 1 and is strictly less than 1 at $\alpha = 1/2$ unless $P_0 = Q$. In such a case there appears to be no loss in generality to work with the choice $\alpha = 1/2$ only, leading to the distance $d$ as in Lemma 2.3. This lemma shows that this situation arises if $\mathscr{P}$ is convex. □

The following theorem translates Theorem 6.2 into the form needed for the proof of our main results. Recall the definition of the covering numbers for testing $N_t(\varepsilon, \mathscr{P}, d; P_0, P^*)$ in (2.2).

THEOREM 6.3. *Suppose $P^* \in \mathscr{P}$ and $P_0(p/p^*) < \infty$ for all $P \in \mathscr{P}$. Assume that there exists a nonincreasing function $D$ such that, for some $\varepsilon_n \geq 0$ and every $\varepsilon > \varepsilon_n$,*

(6.5) $$N_t(\varepsilon, \mathscr{P}, d; P_0, P^*) \leq D(\varepsilon).$$

*Then for every $\varepsilon > \varepsilon_n$ there exists a test $\phi_n$ (depending on $\varepsilon > 0$) such that, for every $J \in \mathbb{N}$,*

(6.6)
$$P_0^n \phi_n \leq D(\varepsilon) \frac{e^{-n\varepsilon^2/4}}{1 - e^{-n\varepsilon^2/4}},$$
$$\sup_{\{P \in \mathscr{P} : d(P, P^*) > J\varepsilon\}} Q(P)^n (1 - \phi_n) \leq e^{-nJ^2\varepsilon^2/4}.$$

PROOF. Define $\mathscr{Q}$ as the set of all finite measures $Q(P)$ as $P$ ranges over $\mathscr{P}$ (where $p_0/p^* = 0$ if $p_0 = 0$) and define $\tau(Q_1, Q_2) = d(P_1, P_2)$. Then $Q(P^*) = P_0$ and, hence, $d(P, P^*) = \tau(Q(P), P_0)$. Identify $P$ of Theorem 6.2 with the present measure $P_0$. By the definitions (2.2) and (6.2), we have $N_\tau(\varepsilon, \mathscr{Q}) \leq N_t(\varepsilon, \mathscr{P}, d; P_0, P^*) \leq D(\varepsilon)$ for every $\varepsilon > \varepsilon_n$. Therefore, the test function guaranteed to exist by Theorem 6.2 satisfies

$$P_0^n \phi_n \leq \sum_{j=1}^\infty D(j\varepsilon) e^{-nj^2\varepsilon^2/4} \leq D(\varepsilon) \sum_{j=1}^\infty e^{-nj^2\varepsilon^2/4},$$

because $D$ is nonincreasing. This can be bounded further (as in the assertion) since, for all $0 < x < 1$, $\sum_{n\geq 1} x^{n^2} \leq x/(1-x)$. The second line in the assertion is simply the second line in (6.3). $\square$

**7. Proofs of the main theorems.** The following lemma is analogous to Lemma 8.1 in [8] and can be proved in the same manner.

LEMMA 7.1. *For given $\varepsilon > 0$ and $P^* \in \mathscr{P}$ such that $P_0 \log(p_0/p^*) < \infty$, define $B(\varepsilon, P^*; P_0)$ by (2.3). Then for every $C > 0$ and probability measure $\Pi$ on $\mathscr{P}$,*

$$P_0^n\left(\int \prod_{i=1}^n \frac{p}{p^*}(X_i)\, d\Pi(P) < \Pi(B(\varepsilon, P^*; P_0)) e^{-n\varepsilon^2(1+C)}\right) \leq \frac{1}{C^2 n\varepsilon^2}.$$

PROOF OF THEOREM 2.2. In view of (2.5), the conditions of Theorem 6.3 are satisfied, with the function $D(\varepsilon) = e^{n\varepsilon_n^2}$ (i.e., constant in $\varepsilon > \varepsilon_n$). Let $\phi_n$ be the test as in the assertion of this theorem for $\varepsilon = M\varepsilon_n$ and $M$ a large constant, to be determined later in the proof.

For $C > 0$, also to be determined later in the proof, let $\Omega_n$ be the event

$$(7.1) \qquad \int_{\mathscr{P}} \prod_{i=1}^n \frac{p}{p^*}(X_i)\, d\Pi(P) \geq e^{-(1+C)n\varepsilon_n^2} \Pi(B(\varepsilon_n, P^*; P_0)).$$

Then $P_0^n(\mathscr{X}^n \setminus \Omega_n) \leq 1/(C^2 n\varepsilon_n^2)$, by Lemma 7.1.

Set $\hat{\Pi}_n(\varepsilon) = \Pi_n(P \in \mathscr{P} : d(P, P^*) > \varepsilon | X_1, \ldots, X_n)$. For every $n \geq 1$ and $J \in \mathbb{N}$, we can decompose

$$(7.2) \quad \begin{aligned} P_0^n \hat{\Pi}_n(JM\varepsilon_n) &= P_0^n(\hat{\Pi}_n(JM\varepsilon_n)\phi_n) + P_0^n(\hat{\Pi}_n(JM\varepsilon_n)(1-\phi_n)\mathbb{1}_{\Omega_n^c}) \\ &\quad + P_0^n(\hat{\Pi}_n(JM\varepsilon_n)(1-\phi_n)\mathbb{1}_{\Omega_n}). \end{aligned}$$

We estimate the three terms on the right-hand side separately. Because $\hat{\Pi}_n(\varepsilon) \leq 1$, the middle term is bounded by $1/(C^2 n\varepsilon_n^2)$. This converges to zero as $n\varepsilon_n^2 \to \infty$ for fixed $C$ and/or can be made arbitrarily small by choosing a large constant $C$ if $n\varepsilon_n^2$ is bounded away from zero.

By the first inequality in (6.6), the first term on the right-hand side of (7.2) is bounded by

$$P_0^n(\hat{\Pi}_n(JM\varepsilon_n)\phi_n) \leq P_0^n \phi_n \leq \frac{e^{(1-M^2/4)n\varepsilon_n^2}}{1 - e^{-nM^2\varepsilon_n^2/4}}.$$

For sufficiently large $M$, the expression on the right-hand side is bounded above by $2e^{-n\varepsilon_n^2 M^2/8}$ for sufficiently large $n$ and, hence, can be made arbitrarily small by choice of $M$, or converges to 0 for fixed $M$ if $n\varepsilon_n^2 \to \infty$.

Estimation of the third term on the right-hand side of (7.2) is more involved. Because $P_0(p^* > 0) = 1$, we can write

$$
\begin{aligned}
(7.3) \quad & P_0^n(\hat{\Pi}_n(JM\varepsilon_n)(1-\phi_n)\mathbb{1}_{\Omega_n}) \\
& = P_0^n(1-\phi_n)\mathbb{1}_{\Omega_n}\left[\frac{\int_{d(P,P^*)>JM\varepsilon_n}\prod_{i=1}^n (p/p^*)(X_i)\,d\Pi(P)}{\int_{\mathscr{P}}\prod_{i=1}^n (p/p^*)(X_i)\,d\Pi(P)}\right],
\end{aligned}
$$

where we have written the arguments $X_i$ for clarity. By the definition of $\Omega_n$, the integral in the denominator is bounded below by $e^{-(1+C)n\varepsilon_n^2}\Pi(B(\varepsilon_n, P^*; P_0))$. Inserting this bound, writing $Q(P)$ for the measure defined by $dQ(P) = (p_0/p^*)\,dP$, and using Fubini's theorem, we can bound the right-hand side of the preceding display by

$$
(7.4) \quad \frac{e^{(1+C)n\varepsilon_n^2}}{\Pi(B(\varepsilon_n, P^*; P_0))}\int_{d(P,P^*)>JM\varepsilon_n} Q(P)^n(1-\phi_n)\,d\Pi(P).
$$

Setting $\mathscr{P}_{n,j} = \{P \in \mathscr{P}: M\varepsilon_n j < d(P, P^*) \leq M\varepsilon_n(j+1)\}$, we can decompose $\{P: d(P, P^*) > JM\varepsilon_n\} = \bigcup_{j \geq J} \mathscr{P}_{n,j}$. The tests $\phi_n$ have been chosen to satisfy the inequality $Q(P)^n(1-\phi_n) \leq e^{-nj^2M^2\varepsilon_n^2/4}$ uniformly in $P \in \mathscr{P}_{n,j}$. [Cf. the second inequality in (6.6).] It follows that the preceding display is bounded by

$$
\begin{aligned}
& \frac{e^{(1+C)n\varepsilon_n^2}}{\Pi(B(\varepsilon_n, P^*; P_0))}\sum_{j\geq J} e^{-nj^2M^2\varepsilon_n^2/4}\Pi(\mathscr{P}_{n,j}) \\
& \leq \sum_{j\geq J} e^{(1+C)n\varepsilon_n^2 + n\varepsilon_n^2 M^2 j^2/8 - nj^2 M^2\varepsilon_n^2/4},
\end{aligned}
$$

by (2.11). For fixed $C$ and sufficiently large $M$, this converges to zero if $n\varepsilon_n^2$ is bounded away from zero and $J = J_n \to \infty$. □

PROOF OF THEOREM 2.1. Because $\Pi$ is a probability measure, the numerator in (2.11) is bounded above by 1. Therefore, the prior mass condition (2.11) is implied (for large $j$) by the prior mass condition (2.4). We conclude that the assertion of Theorem 2.1, but with $M = M_n \to \infty$, follows from Theorem 2.2. That in fact it suffices that $M$ is sufficiently large follows by inspection of the preceding proof. □

PROOF OF THEOREM 2.4. The proof of this theorem follows the same steps as the preceding proofs. A difference is that we cannot appeal to the preparatory lemmas and theorems to split the proof in separate steps. The shells $\mathscr{P}_{n,j} = \{P \in P: Mj\varepsilon_n < d(P, \mathscr{P}^*) < M(j+1)\varepsilon_n\}$ must be covered by sets $B_{n,j,i}$ as in the definition (2.16), and for each such set we use the appropriate element $P_{n,j,i}^* \in \mathscr{P}^*$ to define a test $\phi_{n,j,i}$ and to rewrite the left-hand side of (7.3). We omit the details. □

**8. Technical lemmas.** Lemma 8.1 is used to upper bound the Kullback–Leibler divergence and the expectation of the squared logarithm by a function of the $L_1$-norm. A similar lemma was presented in [16], where both $p$ and $q$ were assumed to be densities of probability distributions. We generalize this result to the case where $q$ is a finite measure and we are forced to use the $L_1$ instead of the Hellinger distance.

LEMMA 8.1. *For every $b > 0$, there exists a constant $\varepsilon_b > 0$ such that, for every probability measure $P$ and finite measure $Q$ with $0 < h^2(p,q) < \varepsilon_b P(p/q)^b$,*

$$P \log \frac{p}{q} \lesssim h^2(p,q)\left(1 + \frac{1}{b}\log_+ \frac{1}{h(p,q)} + \frac{1}{b}\log_+ P\left(\frac{p}{q}\right)^b\right) + \|p-q\|_1,$$

$$P\left(\log \frac{p}{q}\right)^2 \lesssim h^2(p,q)\left(1 + \frac{1}{b}\log_+ \frac{1}{h(p,q)} + \frac{1}{b}\log_+ P\left(\frac{p}{q}\right)^b\right)^2.$$

PROOF. The function $r:(0,\infty) \to \mathbb{R}$ defined implicitly by $\log x = 2(\sqrt{x}-1) - r(x)(\sqrt{x}-1)^2$ possesses the following properties:

- $r$ is nonnegative and decreasing.
- $r(x) \sim \log(1/x)$ as $x \downarrow 0$, whence there exists $\varepsilon' > 0$ such that $r(x) \leq 2\log(1/x)$ on $[0,\varepsilon']$. (A computer graph indicates that $\varepsilon' = 0.4$ will do.)
- For every $b > 0$, there exists $\varepsilon_b'' > 0$ such that $x^b r(x)$ is increasing on $[0, \varepsilon_b'']$. (For $b \geq 1$, we may take $\varepsilon_b'' = 1$, but for $b$ close to zero, $\varepsilon_b''$ must be very small.)

In view of the definition of $r$ and the first property, we can write

$$P \log \frac{p}{q} = -2P\left(\sqrt{\frac{q}{p}} - 1\right) + \Pr\left(\frac{q}{p}\right)\left(\sqrt{\frac{q}{p}} - 1\right)^2$$

$$\leq h^2(p,q) + 1 - \int q \, d\mu + \Pr\left(\frac{q}{p}\right)\left(\sqrt{\frac{q}{p}} - 1\right)^2$$

$$\leq h^2(p,q) + \|p-q\|_1 + r(\varepsilon)h^2(p,q) + \Pr\left(\frac{q}{p}\right)\mathbb{1}\left\{\frac{q}{p} \leq \varepsilon\right\},$$

for any $0 < \varepsilon \leq 4$, where we use the fact that $|\sqrt{q/p} - 1| \leq 1$ if $q/p \leq 4$. Next we choose $\varepsilon \leq \varepsilon_b''$ and use the third property to bound the last term on the right-hand side by $P(p/q)^b \varepsilon^b r(\varepsilon)$. Combining the resulting bound with the second property, we then obtain, for $\varepsilon \leq \varepsilon' \wedge \varepsilon_b'' \wedge 4$,

$$P \log \frac{p}{q} \leq h^2(p,q) + \|p-q\|_1 + 2\log \frac{1}{\varepsilon} h^2(p,q) + 2\varepsilon^b \log \frac{1}{\varepsilon} P\left(\frac{p}{q}\right)^b.$$

For $\varepsilon^b = h^2(p,q)/P(p/q)^b$, the second and third terms on the right-hand side take the same form. If $h^2(p,q) < \varepsilon_b P(p/q)^b$ for a sufficiently small $\varepsilon_b$, then this choice is eligible and the first inequality of the lemma follows. Specifically, we can choose $\varepsilon_b \leq (\varepsilon' \wedge \varepsilon_b'' \wedge 4)^b$.

To prove the second inequality, we first note that, since $|\log x| \leq 2|\sqrt{x} - 1|$ for $x \geq 1$,

$$P\left(\log\frac{p}{q}\right)^2 \mathbb{1}\left\{\frac{q}{p} \geq 1\right\} \leq 4P\left(\sqrt{\frac{q}{p}} - 1\right)^2 = 4h^2(p,q).$$

Next, with $r$ as in the first part of the proof,

$$P\left(\log\frac{p}{q}\right)^2 \mathbb{1}\left\{\frac{q}{p} \leq 1\right\} \leq 8P\left(\sqrt{\frac{q}{p}} - 1\right)^2 + 2Pr^2\left(\frac{q}{p}\right)\left(\sqrt{\frac{q}{p}} - 1\right)^4 \mathbb{1}\left\{\frac{q}{p} \leq 1\right\}$$

$$\leq 8h^2(p,q) + 2r^2(\varepsilon)h^2(p,q) + 2\varepsilon^b r^2(\varepsilon)P\left(\frac{p}{q}\right)^b,$$

for $\varepsilon \leq \varepsilon_{b/2}''$, in view of the third property of $r$. (The power of 4 in the first line of the array can be lowered to 2 or 0, as $|\sqrt{q/p} - 1| \leq 1$.) We can use the second property of $r$ to bound $r(\varepsilon)$ and next choose $\varepsilon^b = h^2(p,q)/P(p/q)^b$ to finish the proof. Specifically, we can choose $\varepsilon_b \leq (\varepsilon' \wedge \varepsilon_{b/2}'')^b$. □

LEMMA 8.2. *If $p, p_n, p_\infty$ are probability densities in $L_1(\mu)$ such that $p_n \to p_\infty$ as $n \to \infty$, then $\liminf_{n\to\infty} P\log(p/p_n) \geq P\log(p/p_\infty)$.*

PROOF. If $X_n = p_n/p$ and $X = p_\infty/p$, then $X_n \to X$ in $P$-probability and in mean. We can write $P\log(p_n/p)$ as the sum of $P(\log X_n)\mathbb{1}_{X_n > 1}$ and $P(\log X_n)\mathbb{1}_{X_n < 1}$. Because $0 \leq (\log x)\mathbb{1}_{x > 1} \leq x$, the sequence $(\log X_n)\mathbb{1}_{X_n > 1}$ is dominated in absolute value by the sequence $|X_n|$, and, hence, is uniformly integrable. By a suitable version of the dominated convergence theorem, we have $P(\log X_n)\mathbb{1}_{X_n > 1} \to P(\log X_\infty)\mathbb{1}_{X_\infty > 1}$. Because the variables $(\log X_n)\mathbb{1}_{X_n < 1}$ are nonnegative, we can apply Fatou's lemma to see that $\limsup P(\log X_n)\mathbb{1}_{X_n < 1} \leq P(\log X_\infty)\mathbb{1}_{X_\infty \leq 1}$. □

## REFERENCES

[1] BERK, R. H. (1966). Limiting behavior of posterior distributions when the model is incorrect. *Ann. Math. Statist.* **37** 51–58. MR0189176 [Corrigendum **37** 745–746.]
[2] BIRGÉ, L. (1983). Approximation dans les espaces métriques et théorie de l'estimation. *Z. Wahrsch. Verw. Gebiete* **65** 181–238. MR0722129
[3] BUNKE, O. and MILHAUD, X. (1998). Asymptotic behavior of Bayes estimates under possibly incorrect models. *Ann. Statist.* **26** 617–644. MR1626075
[4] DIACONIS, P. and FREEDMAN, D. (1986). On the consistency of Bayes estimates (with discussion). *Ann. Statist.* **14** 1–67. MR0829555

[5] DIACONIS, P. and FREEDMAN, D. (1986). On inconsistent Bayes estimates of location. *Ann. Statist.* **14** 68–87. MR0829556

[6] FERGUSON, T. S. (1973). A Bayesian analysis of some non-parametric problems. *Ann. Statist.* **1** 209–230. MR0350949

[7] FERGUSON, T. S. (1974). Prior distributions on spaces of probability measures. *Ann. Statist.* **2** 615–629. MR0438568

[8] GHOSAL, S., GHOSH, J. K. and VAN DER VAART, A. W. (2000). Convergence rates of posterior distributions. *Ann. Statist.* **28** 500–531. MR1790007

[9] GHOSAL, S. and VAN DER VAART, A. W. (2001). Entropies and rates of convergence for maximum likelihood and Bayes estimation for mixtures of normal densities. *Ann. Statist.* **29** 1233–1263. MR1873329

[10] LE CAM, L. M. (1986). *Asymptotic Methods in Statistical Decision Theory.* Springer, New York. MR0856411

[11] MEGGINSON, R. E. (1998). *An Introduction to Banach Space Theory.* Springer, New York. MR1650235

[12] PFANZAGL, J. (1988). Consistency of maximum likelihood estimators for certain nonparametric families, in particular: Mixtures. *J. Statist. Plann. Inference* **19** 137–158. MR0944202

[13] SCHWARTZ, L. (1965). On Bayes procedures. *Z. Wahrsch. Verw. Gebiete* **4** 10–26. MR0184378

[14] SHEN, X. and WASSERMAN, L. (2001). Rates of convergence of posterior distributions. *Ann. Statist.* **29** 687–714. MR1865337

[15] STRASSER, H. (1985). *Mathematical Theory of Statistics.* de Gruyter, Berlin. MR0812467

[16] WONG, W. H. and SHEN, X. (1995). Probability inequalities for likelihood ratios and convergence rates of sieve MLE's. *Ann. Statist.* **23** 339–362. MR1332570

DEPARTMENT OF MATHEMATICS
VRIJE UNIVERSITEIT AMSTERDAM
DE BOELELAAN 1081A
1081 HV AMSTERDAM
THE NETHERLANDS
E-MAIL: aad@cs.vu.nl
kleijn@cs.vu.nl